\documentclass[12pt]{amsart}
\usepackage{amsbsy}
\usepackage{graphicx,epsfig,subfigure,psfrag}
\begin{document}

\newtheorem{theorem}{Theorem}
\newtheorem{proposition}{Proposition}
\newtheorem{lemma}{Lemma}
\newtheorem{corollary}{Corollary}
\newtheorem{definition}{Definition}
\newtheorem{remark}{Remark}
\newcommand{\tex}{\textstyle}
\numberwithin{equation}{section} \numberwithin{theorem}{section}
\numberwithin{proposition}{section} \numberwithin{lemma}{section}
\numberwithin{corollary}{section}
\numberwithin{definition}{section} \numberwithin{remark}{section}
\newcommand{\ren}{\mathbb{R}^N}
\newcommand{\re}{\mathbb{R}}
\newcommand{\n}{\nabla}
\newcommand{\iy}{\infty}
\newcommand{\pa}{\partial}
\newcommand{\fp}{\noindent}
\newcommand{\ms}{\medskip\vskip-.1cm}
\newcommand{\mpb}{\medskip}
\newcommand{\AAA}{{\bf A}}
\newcommand{\BB}{{\bf B}}
\newcommand{\CC}{{\bf C}}
\newcommand{\DD}{{\bf D}}
\newcommand{\EE}{{\bf E}}
\newcommand{\FF}{{\bf F}}
\newcommand{\GG}{{\bf G}}
\newcommand{\oo}{{\mathbf \omega}}
\newcommand{\Am}{{\bf A}_{2m}}
\newcommand{\CCC}{{\mathbf  C}}
\newcommand{\II}{{\mathrm{Im}}\,}
\newcommand{\RR}{{\mathrm{Re}}\,}
\newcommand{\eee}{{\mathrm  e}}
\newcommand{\LL}{L^2_\rho(\ren)}
\newcommand{\LLL}{L^2_{\rho^*}(\ren)}
\renewcommand{\a}{\alpha}
\renewcommand{\b}{\beta}
\newcommand{\g}{\gamma}
\newcommand{\G}{\Gamma}
\renewcommand{\d}{\delta}
\newcommand{\D}{\Delta}
\newcommand{\e}{\varepsilon}
\newcommand{\var}{\varphi}
\newcommand{\lll}{\l}
\renewcommand{\l}{\lambda}
\renewcommand{\o}{\omega}
\renewcommand{\O}{\Omega}
\newcommand{\s}{\sigma}
\renewcommand{\t}{\tau}
\renewcommand{\th}{\theta}
\newcommand{\z}{\zeta}
\newcommand{\wx}{\widetilde x}
\newcommand{\wt}{\widetilde t}
\newcommand{\noi}{\noindent}
\newcommand{\uu}{{\bf u}}
\newcommand{\xx}{{\bf x}}
\newcommand{\yy}{{\bf y}}
\newcommand{\zz}{{\bf z}}
\newcommand{\aaa}{{\bf a}}
\newcommand{\cc}{{\bf c}}
\newcommand{\jj}{{\bf j}}
\newcommand{\ggg}{{\bf g}}
\newcommand{\UU}{{\bf U}}
\newcommand{\YY}{{\bf Y}}
\newcommand{\HH}{{\bf H}}
\newcommand{\GGG}{{\bf G}}
\newcommand{\VV}{{\bf V}}
\newcommand{\ww}{{\bf w}}
\newcommand{\vv}{{\bf v}}
\newcommand{\hh}{{\bf h}}
\newcommand{\di}{{\rm div}\,}
\newcommand{\ii}{{\rm i}\,}
\newcommand{\inA}{\quad \mbox{in} \quad \ren \times \re_+}
\newcommand{\inB}{\quad \mbox{in} \quad}
\newcommand{\inC}{\quad \mbox{in} \quad \re \times \re_+}
\newcommand{\inD}{\quad \mbox{in} \quad \re}
\newcommand{\forA}{\quad \mbox{for} \quad}
\newcommand{\whereA}{,\quad \mbox{where} \quad}
\newcommand{\asA}{\quad \mbox{as} \quad}
\newcommand{\andA}{\quad \mbox{and} \quad}
\newcommand{\withA}{,\quad \mbox{with} \quad}
\newcommand{\orA}{,\quad \mbox{or} \quad}
\newcommand{\atA}{\quad \mbox{at} \quad}
\newcommand{\onA}{\quad \mbox{on} \quad}
\newcommand{\ef}{\eqref}
\newcommand{\ssk}{\smallskip}
\newcommand{\LongA}{\quad \Longrightarrow \quad}
\def\com#1{\fbox{\parbox{6in}{\texttt{#1}}}}
\def\N{{\mathbb N}}
\def\A{{\cal A}}
\newcommand{\de}{\,d}
\newcommand{\eps}{\varepsilon}
\newcommand{\be}{\begin{equation}}
\newcommand{\ee}{\end{equation}}
\newcommand{\spt}{{\mbox spt}}
\newcommand{\ind}{{\mbox ind}}
\newcommand{\supp}{{\mbox supp}}
\newcommand{\dip}{\displaystyle}
\newcommand{\prt}{\partial}
\renewcommand{\theequation}{\thesection.\arabic{equation}}
\renewcommand{\baselinestretch}{1.1}
\newcommand{\Dm}{(-\D)^m}

\title
{\bf  Boundary characteristic point regularity for\\semilinear
 reaction-diffusion equations:\\ towards an ODE criterion}

\author {V.A.~Galaktionov and V.~Maz'ya}

\address{Department of Mathematical Sciences, University of Bath,
 Bath BA2 7AY, UK}
\email{vag@maths.bath.ac.uk}

\address{Department of Mathematical Sciences, University of Liverpool,
  M$\&$O Building,
Liverpool, L69 3BX, UK \,\, {\rm and} \,\, Department of
Mathematics, Link\"oping University, SE-58183, Link\"oping,
Sweden}
 \email{vlmaz@liv.ac.uk \, and \,  vlmaz@mai.liu.se}


\date{\today}


\begin{abstract}

The classic  problem of regularity of boundary characteristic
points for
 semilinear  heat equations with homogeneous Dirichlet conditions is considered.
It is shown that famous
 Petrovskii's (the so-called, $2\sqrt{\log\, \log}$) criterion of boundary regularity for the heat equation (1934)
  can be adapted  to  classes of semilinear parabolic equations of reaction-diffusion
  type,
  and now this takes the form of  an {\em ODE regularity criterion}.
   Namely, after a special matching with a {\em Boundary
  Layer}, the regularity problem reduces
  to
 a one-dimensional  perturbed nonlinear dynamical system  for the first Fourier-like coefficient
 of the solution in an {\em Inner Region}.

A similar ODE criterion, with an analogous matching procedures, is
shown formally to exist for semilinear fourth-order bi-harmonic
equations of reaction-diffusion type.
 Extensions to regularity problems of backward paraboloid vertexes  in $\ren$ are
discussed.

\end{abstract}

\maketitle





\section{Introduction: main semilinear equations, classic regularity problem, results, and layout}
\label{S0}


\subsection{Semilinear reaction-diffusion PDEs near parabola vertexes}

The present paper is devoted to a  systematic study of the
regularity of the origin $(0,0)$ as a  boundary point for
semilinear heat (reaction-diffusion) equations, which we first
consider in 1D:
 \be
 \label{sem1}
 u_t= u_{xx} + f(x,t,u) \inB Q_0 \subset \re \times [-1,0), \quad u=0 \,\,\, \mbox{on} \,\,\, \pa Q_0,
 \quad u(x,-1)=u_0(x).
  \ee
  Here $Q_0$ is a sufficiently smooth domain  such that  $(0,0) \in
  \overline{\partial Q_0}$ ($\pa Q_0$ denotes the lateral boundary of $Q_0$)
    is its only {\em characteristic} boundary point, i.e., in the $\{x,t\}$-plane,

    (i) the straight line $\{t=0\}$
is tangent to $\partial Q_0$ at this point, and

(ii) no such points exist on $\pa Q_0$ for $t \in [-1,0)$.

\noi According to \ef{sem1}, we pose the zero Dirichlet condition
on the lateral boundary $\pa Q_0$ and prescribe arbitrary bounded
initial data $u_0(x)$ at $t=-1$ in $Q_0 \cap \{t=-1\}$.

We naturally assume that nonlinearities  $f(x,t,u)$ in \ef{sem1}
satisfy necessary regularity and growth in $u$ hypotheses that
guarantee existence and uniqueness of a smooth classical solution
$u \in C^{2,1}_{x,t}(Q_0) \cap C(\bar Q_0)$\footnote{As customary
in parabolic PDEs \cite{Fr0}, the closure $\bar Q_0$ does not
include the ``upper lid", which is the vertex $(0,0)$ only.} of
\ef{sem1} in $Q_0$ by classic parabolic theory; see e.g.,
well-known monographs \cite{A1,EidSys, Fr0, Fr}. A standard
regularity issue in general PDE theory is then as follows: how to
control the solution at the vertex, i.e., the problem is:
 $$
 u(0,0^-)=?
 $$
 Note that this does not {\em a priori} exclude {\em blow-up} at
  the vertex (regardless zero Dirichlet conditions on the lateral boundary) where
  $
  |u(0,0^-)|=+ \iy$, in the sense of $\limsup$.

 More precisely, as customary in regularity theory, the goal is
to derive conditions {\em showing how given smooth nonlinear
perturbations $f(\cdot)$ can affect regularity of the vertex
$(0,0)$ of a such a {backward parabola} $\pa Q_0$}. The regularity
of $(0,0)$
 (in Wiener's sense)
 means:
 \be
 \label{Reg11}
 {\bf regularity\,\, of\,\, (0,0)}: \quad u(0,0^-)=0 \quad \mbox{for
any initial data $u_0$}.
 \ee



\ssk

 Recall that, as is well known, for nonlinearities
  $f(\cdot) \equiv 0$, i.e., for the pure  {\em heat equation}
 \be
 \label{heat1}
 u_t=  u_{xx} \inB Q_0,  \quad u=0 \,\,\, \mbox{on} \,\,\, \pa Q_0,\quad u(x,-1)=u_0(x),
  \ee
 this regularity problem was solved by Petrovskii in 1934
 \cite{Pet34Cr, Pet35}, who introduced his famous
{\em Petrovskii's regularity criterion} (the so-called
``$2\sqrt{\log \log}$--one"; see details below).

\ssk

 Indeed, many and often strong and delicate boundary regularity and related
 asymptotic
results are now known for a number of quasilinear parabolic
equations,
including even a few for  degenerate porous medium operators.
Nevertheless, some difficult questions remain open even for the
second-order parabolic equations with order-preserving semigroups.
 We refer to the results and surveys in \cite{Abd00, Abd03, Abd05, GalPet2m, Herr04, Mam01, Moss00}
as a guide to a full history and already existing interesting
extensions of these important results. Concerning further
developing of Wiener's ideas in linear parabolic equations, see
references and results in \cite{Land69, Lanc73, Evans82} and in
\cite{Wat97}. However, a more systematic study of those regularity
issues for equations such as \ef{sem1} with rather general
nonlinear perturbations $f(\cdot)$ was not done properly still. In
fact,  it turned out that, for such arbitrary $f(\cdot)$'s,
 the classical barrier methods hardly applied and another asymptotic approach
was necessary. We propose this in the present paper for a wide
class of semilinear parabolic PDEs.

\subsection{Layout of the paper: key models, nonlinearities, and extensions}

Section \ref{S1}  contain some preliminary discussions and
results. In Sections \ref{S3}--\ref{S3.4}, the main goal is to
show how a general ``nonlinear perturbation" $f(\cdot)$ in
\ef{sem1} affects the regularity conditions by deriving sharp
asymptotics of solutions near characteristic points. To this end,
we apply a method of a matched asymptotic (blow-up) expansion,
where the {\em Boundary Layer} behaviour close to the lateral
boundary $\partial Q_0$ (Section \ref{S3.3}) is matched, as $t \to
0^-$, with a {\em centre subspace behaviour} in an {\em Inner
Region} (Section \ref{S3.4}). This leads to a nonlinear dynamical
system for the first Fourier coefficient in the eigenfunction
expansion via standard Hermite polynomials as eigenfunctions of
the linear Hermite operator obtained after blow-up scaling near
the vertex. Overall, the vertex regularity is shown to be governed
by an {\em ODE criterion}, which principally does not admit any
simply integral (Osgood--Dini-type) treatment as in Petrovskii's
one.

 Indeed, such an approach falls into the scope of
typical ideas of asymptotic PDE theory, which got a full
mathematical justification for many problems of interest. In
particular, we refer to a recent general asymptotic analysis
performed in \cite{KozMaz04}. According to its classification, our
matched blow-up approach corresponds to perturbed {\sc
one-dimensional} dynamical systems, i.e., to a rather elementary
case being however a constructive one that detects a number of new
asymptotic/regularity results.

In particular, to show a typical ``interaction" between the linear
Laplacian and the nonlinear perturbations in \ef{sem1}, we
initially concentrate on the simplest case, with
 \be
 \label{3.1}
  \tex{
 f(\cdot)= \frac 1{(-t)} \, \kappa(u) u \forA t \in [-1,0),
 }
  \ee
   where $\kappa(u)$ is a smooth enough  function satisfying
   \be
   \label{kap1}
   \kappa(u) \to 0\,\, \mbox{as}\,\, u \to 0,\quad |\kappa(u)| \le 1,  \quad
    \kappa(u) \not = 0 \forA u \ne 0.
 \ee
We show that the nonlinear perturbation \ef{3.1} will then affect
Petrovskii's $2\sqrt{\log \, \log}$-time-factor starting from some
awkward looking  functions such as
 \be
 \label{kap99}
  \tex{
   \kappa(u) \sim |\ln |u| \,|^{\frac 13} {\mathrm e}^{-(3 \sqrt{\pi}|\ln
   |u| \,|)^{2/3}}\to 0 \asA u \to 0.
   }
   \ee
 For more general nonlinearities, we derive a so-called {\em ODE regularity criterion}
of the vertex $(0,0)$, meaning that a special nonlinear ODE for
the first Fourier coefficients of rescaled solutions takes
responsibility for the vertex regularity/irregularity.

 The present research has
been inspired by the regularity study of {\em quasilinear elliptic
equations} with quadratic gradient-dependent nonlinearities
\cite{MS04}, where, in 2D, new asymptotics of solutions near
corner points were discovered.
We also refer to monographs \cite{Gris85, KMR1, KozMaz01, Maz2010,
MazSol2010} and \cite{KozMaz05}, \cite{MayMaz09}--\cite{Maz02} as
an update guide to elliptic regularity theory including
higher-order equations. Sharp asymptotics of solutions of the heat
equation in domains with conical points were derived in
\cite{Gris92, KozMazHE1, KozMazHE2, Koz88}. Higher-order parabolic
equations were treated in \cite{Koz2mSib, Koz2m88}.
 It turned out that, unlike the
present study, such asymptotics are of a self-similar form.
 See also \cite{Kweon07} for a good short survey including
  compressible/incompressible
 Stokes and Navier--Stokes problems.

\ssk

 Therefore, as a next key model
regularity problem,
 we briefly reflect the main differences and
difficulties, which occur by studying the regularity issues for
parabolic equations with a typical {\em quadratic gradient
dependence} in the nonlinear term:
 \be
 \label{quad1}
 u_t=u_{xx} + \kappa(u)u \,(u_x)^2 \inB Q_0,  \quad u=0 \,\,\, \mbox{on} \,\,\, \pa Q_0, \quad u(x,-1)=u_0(x).
  \ee
Then the ODE regularity criterion is expressed in terms of another
1D dynamical system, with a weaker nonlinearity.
We then convincingly show that, for any $\kappa(u)$ in \ef{quad1},
satisfying \ef{kap1}, Petrovskii's linear regularity criterion
takes place, i.e., remains the same as for the heat equation
\ef{heat1}.

\ssk

We also   pay some attention to extensions to similar regularity
problems in domains $Q_0 \subset \ren \times [-1,1)$, with $\pa
Q_0$ having a {\em backward paraboloid shape} and the vertex
$(0,0)$ being their characteristic point. In Section \ref{S3}, we
thus discuss the semilinear problems:
 \be
 \label{sem1N}
 u_t= \D u +
 \left\{
 \begin{matrix}
 \frac{\kappa(u) u}{(-t)}\ssk\\
 \kappa(u)u |\n u|^2
 \end{matrix}
 \right.
  \inB Q_0, \quad u=0 \,\,\, \mbox{on} \,\,\, \pa Q_0.
  \ee

\ssk

Finally, in Appendix B (Appendix A is devoted to the corresponding
spectral theory of rescaled operators), we show how our approach
can be extended to higher-order PDEs, e.g., for the {\em
semilinear bi-harmonic equations} having similar nonlinearities,
with also zero Dirichlet conditions on $\partial Q_0$ and bounded
initial data $u_0$ in $Q_0 \cap \{t=-1\}$. The mathematical
analysis becomes much more difficult and we do not justify
rigorously all its main steps such as the boundary layer and
matching with the Inner Region asymptotics. Moreover, the 1D
dynamical system for the first Fourier coefficients becomes also
more delicate and does not admit such a complete analysis, though
some definite conclusions are possible. We must admit that this
part of our study is formal, though some steps are expected to
admit a full justification, which nevertheless can be rather
time-consuming.


\section{
 Petrovskii's $2\sqrt{\log \log}$--criterion of 1934 and some extensions}
 \label{S1}

We need to explain some details of Petrovskii's classic regularity
analysis for the heat equation performed in 1934-35. Following his
study, we consider the one-dimensional case $N=1$, where the
analysis becomes more clear. Moreover, our further extensions
 to bi-harmonic operators (Appendix B)
  will be also performed  for $N=1$, in view
of rather complicated asymptotics occurred, so we are not
interested in involving extra technicalities.

\ssk


After Wiener's pioneering  regularity criterion for the {\em
Laplace equation} in 1924 \cite{Wien24}, I.G.~Petrovskii
\cite{Pet34Cr, Pet35} was the first who
   completed the study
  of
the regularity question for the 1D and 2D {\em heat equation} in a
non-cylindrical domain.
 We formulate his result in a
{\em blow-up manner}, which in fact was already  used by
Petrovskii in 1934 \cite{Pet34Cr}.

He considered the question on an {\em irregular} or {\em regular}
vertex  $(x,t)=(0,0)$
in the initial-boundary value problem (IBVP)
 \be
 \label{ir1}
  \left\{
 \begin{matrix}
   u_t=u_{xx} \quad \mbox{in}
 \quad Q_0=\{|x|<R(t), \,\,\, -1 < t<0\}, \quad R(t) \to 0^+ \,\,\mbox{as}\,\, t \to
 0^-, \ssk\ssk\ssk\\
  \mbox{with bounded smooth data $u(x,0)=u_0(x)$ on $[-R(-1),R(-1)]$}.
  \quad\,\quad\,\,\,\,
   \end{matrix}
    \right.
  \ee
Here the lateral boundary $\{x= \pm R(t), \,\, t \in [-1,0)\}$ is
given by a function $R(t)$ that is  assumed to be positive,
strictly monotone, $C^1$-smooth for all $-1 \le t  <0$ (with
$R'(t) > 0$), and is allowed to have  a singularity of $R'(t)$ at
$t=0^-$ only. The regularity analysis then detects  the value of
$u(x,t)$ at the end ``blow-up" {\em characteristic} point
$(0,0^-)$, to which the domain $Q_0$ ``shrinks" as $t \to 0^-$.

\ssk

\noi{\bf Remark: on first parabolic regularity results for $m=1$
and $m \ge 2$.} It is well-known that,
 for the heat equation, the first existence  of a classical solution (i.e., continuous at
 $(0,0)$)
   was obtained  by
 Gevrey in 1913--14 \cite{Gev13}
 (see Petrovskii's references  in \cite[p.~55]{Pet34Cr} and \cite[p.~425]{Pet35}),
 which assumed that the H\"older exponent of
 $R(t)$ is larger than $\frac 12$. In our setting, at $t=0^-$, this comprises all types
 of boundaries given by the functions:
  \be
  \label{Gev1}
  R(t)=
(-t)^\nu \quad \mbox{with any $\nu > \frac 12$ are regular}
 \quad(\mbox{Gevrey, 1913--14}).
 \ee
For $2m$th-order parabolic poly-harmonic  equations such as
\ef{Lineq} below, a similar result: 
 \be
  \label{Gev1Mih}
  \mbox{for} \quad
  R(t)=
(-t)^{\frac 1{2m}}, \quad \mbox{the problem is uniquely solvable}
 \ee
 was proved\footnote{However, in
Slobodetskii--Sobolev classes, i.e., continuity at $(0,0)$ is not
understood  in the above Wiener classic sense. In fact, for $m=2$,
Wiener's one \ef{u0} fails for the parabola  with $R(t)=
5(-t)^{\frac 1{4}}$ in \ef{Gev1Mih}, while the one with  $R(t)=
4(-t)^{\frac 1{4}}$ remains regular, \cite[\S~4.3]{GalPet2m}.} by
Miha$\breve{\rm i}$lov \cite{Mih61} almost sixty years later and
fifty years ago.

\ssk

\noi{\bf Definitions.}\, \underline{\bf (i) Regular point}: as
usual in potential theory, {\em the point $(x,t)=(0,0)$ is called
{\em regular}} (in Wiener's sense, see \cite{Maz99}),
{\em if any value of the solution $u(x,t)$ can be prescribed there
by continuity as a standard boundary value on $\partial Q_0$}. In
particular,
 as a convenient and key for us  evolution illustration, {\em $(0,0)$ is regular if the
 continuity holds for any initial data $u_0(x)$ in the following sense}:
  \be
  \label{u0}
  u=0 \,\,\, \mbox{at the lateral boundary} \,\,\,\{|x|=R(t), \, -1 \le t<0\} \LongA
  u(0,0^-)=0.
  \ee

\noi\underline{\bf (ii) Irregular point}: otherwise, {\em the
point $(0,0)$ is {\em irregular},
      if the value $u(0,0^-)$ is not
fixed by boundary conditions, i.e., $u(0,0) \ne 0$ for some data
$u_0$}, and hence is given by a ``blow-up
evolution" as $t \to 0^-$. Hence, formally, $(0,0)$ does not
belong to the parabolic boundary of $Q_0$.






\ssk

 {\bf Petrovskii's ``${\mathbf 2} \sqrt{\mathbf {log\, log}}$".} Using
novel barriers as upper and lower solutions of \ef{ir1},
 Petrovskii  \cite{Pet34Cr, Pet35} established the following
 ``$2\sqrt{\log\log}$-criterion":
 \be
 \label{PP1}
 \fbox{$
 \begin{matrix}
{\rm (i)}\,R(t)=2 \sqrt{-t} \,\,\sqrt{\ln|\ln (-t)|} \LongA
\mbox{$(0,0)$ is regular}, \andA \qquad\quad \ssk\\ {\rm (ii)}\,
R(t)=2(1+\e)\sqrt{-t}\,\, \sqrt{ \ln|\ln (-t)|}, \,\, \e>0 \LongA
 \mbox{$(0,0)$ is irregular}.
  \end{matrix}
  $}
  \ee

More precisely, he also showed that, for the curve expressed in
terms of
 a
positive function $\rho(h) \to 0^+$ as $h \to 0^+$ ($\rho(h) \sim
\frac 1{|\ln h | }$ is about right) as follows:
 \be
 \label{RR1NN}
  \tex{
  R(t) = 2 \sqrt{-t} \,\,\sqrt{-\ln \rho(-t)},
   }
  \ee
 the {\em sharp regularity criterion} holds (in Petrovskii's original notation):
  \be
  \label{RR1}
   \fbox{$
\int\limits_0 \frac{{\rho(h)} \sqrt{|\ln \rho(h)|}}h \, {\mathrm
  d}h < \,(=)\, +\iy
   \LongA \mbox{$(0,0^-)$ is irregular (regular)}.
 $}
  \ee

 Both converging (irregularity)
 and diverging (regularity) integrals in \ef{RR1} as Dini--Osgood-type
 regularity criteria already appeared in the first Petrovskii paper
 \cite[p.~56]{Pet34Cr} of 1934. Further historical and
mathematical comments concerning Petrovskii's analysis including
earlier \cite{Khin36} (1933) Khinchin's criterion
 in a probability representation can be
 found in a survey in \cite{GalPet2m}.

 \ssk

Petrovskii's integral criterion of the Dini--Osgood type given in
\ef{RR1}
    is true in the $N$-dimensional
radial case with (see \cite{Abd00, Abd03, Abd05} for a more recent
updating)
  \be
  \label{pp1}
\sqrt{|\ln p(h)|} \quad \mbox{replaced by} \quad {|\ln
p(h)|}^{\frac N 2}.
  \ee

\ssk

It is worth mentioning that, as far as we know, \ef{PP1} is the
first clear appearance of the ``magic" $\sqrt{\rm log\,log}$ in
PDE theory, currently associated with the ``blow-up behaviour" of
the domain $Q_0$ and corresponding solutions. Concerning other
classes of nonlinear PDEs generating blow-up $\sqrt{\rm log\,log}$
in other settings, see references in \cite{GalLog}.

\ssk

Thus, since the 1930s, Petrovskii's regularity
$\sqrt{\log\log}$-factor entered parabolic theory and generated
new types of asymptotic blow-up problems, which have been solved
for a wide class of parabolic equations with variable coefficients
as well as for some quasilinear ones. Nevertheless, such
asymptotic problems were very delicate and some of them  of
Petrovskii's type remained open even in the second-order case,
i.e., for \ef{sem1}, to  be solved in the present paper for the
first time.





 \section{Preliminaries of matched asymptotic expansion}
 \label{S3}

\subsection{The basic initial-boundary value problem}

Thus, we consider the semilinear parabolic equation \ef{sem1},
with  a simple, ``basic" nonlinear perturbation, which we take in
the  separable  form \ef{3.1}.
 The eventual ODE regularity criterion   will then also include
    the behaviour of the nonlinear coefficient  $\kappa(u)$
 as $u \to 0$.

 Hence, our {\em basic second-order initial-boundary value problem} (IBVP) takes the form:
 \be
 \label{ir12}
  \left\{
 \begin{matrix}
   u_t=u_{xx} + \frac 1{(-t)} \, \kappa(u) u\quad \mbox{in}
 \quad Q_0=\{|x|<R(t), \,\,\, -1 < t<0\},
 \ssk\ssk\\
  u=0 \atA x= \pm R(t), \quad -1 \le t <0,\quad\,\,\,
  \qquad\qquad\qquad
  \,\,\,\,\,\,\,\,\,
    \ssk\ssk\\
  u(x,0)=u_0(x) \onA
  [-R(-1),R(-1)],\qquad\quad\,\,\,\,\,\,\,\,\qquad\qquad
   \end{matrix}
    \right.
  \ee
where $u_0(x)$ is a  bounded and smooth function, $u_0(\pm
R(-1))=0$. We then apply to \ef{ir12} both Definitions from the
previous Section \ref{S1}.




 \subsection{Slow growing factor $\var(\t)$}

 According  to \ef{PP1}, we need to assume that
 \be
 \label{ph1}
 R(t)= (-t)^{\frac 12}\, \var(\t) \whereA \t= - \ln(-t) \to + \iy
 \asA t \to 0^-.
  \ee
 Here,  $\var(\t)>0$ is a smooth monotone increasing function
 satisfying $\var'(\t)>0$,
 \be
 \label{vv1}
  \tex{
 \var(\t) \to +\iy, \quad \var'(\t) \to 0^+, \andA \frac {\var'(\t)}{\var(\t)} \to 0 \asA \t \to
 +\iy.
  }
  \ee
 Moreover, as a sharper characterization of the above class of
 {\em slow growing functions}, we use the following criterion:
  \be
  \label{vv2}
   \tex{
   \big( \frac{\var(\t)}{\var'(\t)}\big)' \to \iy \asA \t \to +\iy
    \quad (\var'(\t) \not = 0).
    }
    \ee
    This is a typical condition in blow-up analysis, which distinguishes classes of
    exponential (the limit in \ef{vv2} is 0), power-like (a constant $\not =
    0$), and slow-growing functions. See \cite[pp.~390-400]{SGKM},
    where in Lemma 1 on p.~400, extra properties of slow-growing
    functions \ef{vv2} are proved. For instance, one can use
    a comparison of such a $\var(\t)$ with any power:
     \be
     \label{al1}
     \mbox{for any $\a>0$}, \quad \var(\t) \ll \t^\a \andA
     \var'(\t) \ll \t^{\a-1} \forA \t \gg 1.
     \ee
 Such estimates are useful
 in  evaluating
  perturbation terms in the rescaled equations.

Thus, the monotone positive function $\var(\t)$ in \ef{ph1}
  is assumed to determine a sharp behaviour of the
boundary of $Q_0$ near the shrinking point $(0,0)$ to guarantee
its regularity. In Petrovskii's criterion \ef{PP1}, the almost
optimal
 function, satisfying  \ef{vv1}, \ef{vv2}, is
  \be
  \label{ph2}
   \var_*(\t) = 2\sqrt {\ln \t} \asA \t \to + \iy.
   \ee



\subsection{First kernel scaling and two region expansion}

 By \ef{ph1}, we
perform the similarity scaling
 \be
 \label{ph3}
  \tex{
  u(x,t)=  v(y,\t) \whereA y=\frac x{(-t)^{
  1/2}}.
  }
   \ee
   Then the rescaled function $v(y,\t)$ now solves the
   rescaled IBVP
 \be
 \label{ph4}
  \left\{
  \begin{matrix}
  v_\t= \BB^* v + \kappa(v)v
   \equiv  v_{yy}- \frac 12 \, y v_y + \kappa(v)v
  \inB Q_0=\{|y| < \var(\t), \,\,\, \t>0\}, \ssk\ssk\\
  v=0 \atA y= \pm \var(\t), \,\,\t \ge 0,
  \qquad\qquad\qquad\qquad\qquad\, \qquad\qquad\qquad\qquad
  \ssk\ssk\\
  v(0,y)=v_0(y) \equiv u_0(y) \onA
  [-R(-1),R(-1)].\qquad\quad\,\,\,\,\,\,\,\,\,\,\,\,
   \qquad\qquad\qquad
   \end{matrix}
   \right.
  \ee
 The rescaled equation in \ef{ph4}, for the first time, shows how
 the classic Hermite operator
  \be
  \label{BB*N}
 \tex{
  \BB^*= D^2_y- \frac 12\, y D_y
  }
  \ee
 occurs after blow-up scaling \ef{ph3}.
 By the divergence \ef{vv1}  of $\var(\t) \to +\iy$ as $\t \to +\iy$, it follows that sharp asymptotics of solutions
  will
 essentially depend on the spectral properties of the linear
 operator $\BB^*$ on the whole  line $\re$ (see Appendix A), as well as on the
 nonlinearity $\kappa(v)v$, so that such an asymptotic ``interaction" between linear
 and nonlinear operators therein eventually determines regularity of the vertex.

 Studying asymptotics for the rescaled problem \ef{ph4},
 as usual in asymptotic analysis, this blow-up
 problem is solved by {\em matching of expansions in two regions}:

 \ssk

 (i) In an {\em Inner Region}, which includes arbitrary compact subsets in $y$
 containing the origin $y=0$, and

 \ssk

 (ii) In a {\em Boundary Region} close to the boundaries $y= \pm \var(\t)$, where a {\em  boundary
 layer} occurs.

 \ssk

 \noi Actually, such a two-region structure, with the
 asymptotics specified below, defines the class of generic solutions under
 consideration.
We begin with the simpler analysis in the Boundary Region (ii).

 \section{Boundary layer (BL) theory}
   \label{S3.3}

\subsection{BL-scaling and a perturbed parabolic equation}

 Sufficiently close to the lateral boundary of $Q_0$, it is
 natural to introduce the variables
  \be
  \label{z1}
   \tex{
   z= \frac y{\var(\t)} \andA v(y,\t)=w(z,\t)
    \LongA w_\t= \frac 1{\var^2} \, w_{zz} - \frac 12\, z w_z +
    \frac {\var'}\var \, z w_z + \kappa(w)w.
    }
    \ee
We next introduce the BL-variables
  \be
  \label{z2}
  \tex{
  \xi= \var^{2}(\t)(1-z) \equiv \var(\var-y), \quad \var^{2}(\t) {\mathrm d} \t={\mathrm d}s,
   \andA w(z,\t)= \rho(s) g(\xi,s),
  }
  \ee
  where $\rho(s)>0$ for $s \gg 1$ is an unknown  scaling
   time-factor depending on the
  function $\var(\t)$.
As usual, this $\rho$-scaling is chosen to get uniformly bounded
rescaled solutions, i.e., for nonnegative solutions,
 \be
 \label{g12}
  \tex{
  \sup_{\xi} \,g(\xi,s)=1 \quad \mbox{for all} \quad s \gg 1
  }
  \ee
  (for solutions which remain of changing sign for $s \gg1$, one
  takes $|g(\xi,s)|$ in \ef{g12}).
  By the Strong Maximum Principle (Sturm's Theorem on zero sets,
  see \cite{An1}), $v_y(y,\t)$  has a finite number of zeros in
  $y$ for any $\t >0$ (possible supremum points) and a standard argument
  ensures
  that the normalization \ef{g12} implies that such a $\rho(s)$
  can be treated as sufficiently smooth for $s \gg 1$\footnote{On
  the other hand, one can normalize in \ef{z2} by the smooth function
  $\rho(s) \equiv v(0,\t)$, which also can be regarded as positive
  (negative) for $\t \gg1$ (infinitely many sign changes of $v(0,\t)$ for $\t \gg 1$ also mean that
  $v(y,\t)$ has infinitely many zeros in $y$ that is impossible for the heat equation \cite{An1}).
  This leads  to some slight technical
  differences, though makes the normalization \ef{z5} below more straightforward.}.
This describes  the class of solutions under consideration. For
instance, by the Maximum Principle, it is particular easier to
work out, when:
 \be
 \label{g13}
 (\ref{g12}) \,\,\, \mbox{holds for all nonnegative solutions $u(x,t) \not \equiv 0$}.
 \ee
Respectively, for non-positive solutions, one can use $-1$ as the
normalization  in \ef{g12}.

On substitution into the PDE in \ef{z1}, we obtain the following
small nonlinear perturbation of a linear uniformly parabolic
equation:
 \be
 \label{z3}
  \begin{matrix}
  g_s= \AAA g - \frac 12\,  \frac 1{\var^{2}} \, \xi g_\xi - \frac
  {\var'_\t}{\var}\, \big(1- \frac \xi {\var^{2}}\big) g_\xi \ssk\ssk\\
- 2\, \frac{\var'_\t}{ \var^{3}} \,
  \xi g_\xi -  \frac {\rho'_s}{\rho} \,
  g + \frac 1{\var^2} \,  \kappa(\rho g)\, g, \quad  \mbox{where} \quad  \AAA g=  g'' + \frac 12\, g'.
    \end{matrix}
   \ee
   As usual in boundary layer theory, this means that
we then are looking for a generic
 pattern of the behaviour described by \ef{z3} on compact subsets
 near the lateral boundary,
  \be
  \label{z4}
  |\xi| = o\big(\var^{-2}(\t)\big)
   \LongA |z-1| = o\big(\var^{-4}(\t)\big) \asA \t \to +
   \iy.
   \ee
On these space-time compact subsets, the second term on the
right-hand side of  \ef{z3} becomes asymptotically small, while
all the other linear ones are much smaller in view of the slow
growth/decay assumptions such as \ef{vv2} for $\var(\t)$ and
$\rho(s)$.

\subsection{Passing to the limit and convergence to a BL-profile}

Thus, we arrive at a uniformly parabolic equation \ef{z3}
perturbed by a number of linear and nonlinear terms being, under
given hypothesis, {\em asymptotically small} perturbations of the
stationary elliptic operator $\AAA$. In particular,
 the last nonlinear term in \ef{z3}
 is clearly asymptotically small by the hypotheses \ef{kap1} and
 \ef{vv1}, so that,
for uniformly bounded $g$,
 \be
 \label{nonz1}
  \tex{
   \frac 1{\var^2(\t)} \,\, g \, \kappa(\rho(s) g) \to 0 \asA \t \to
   +\iy.
   }
    \ee



\ssk

The BL representation \ef{z2}, by using the rescaling and
\ef{g12}, naturally leads to
the following asymptotic behaviour at infinity:
 \be
 \label{z5}
  \tex{
   \lim_{s \to +\iy} g(\xi,s) \to 1 \asA \xi \to + \iy,
   }
  \ee
  where all the derivatives also vanish. Then, we arrive at
  the problem of passing to the limit as $s \to + \iy$ in the
  problem \ef{z3}, \ef{z5}. Since, by the definition in \ef{z2},
  the rescaled orbit $\{g(s),\,\, s>0\}$ is uniformly bounded, by
  classic parabolic interior regularity  theory \cite{Fr, EidSys, Eid}, one can pass to the
  limit in \ef{z3} along a subsequence $\{s_k\} \to +\iy$. Namely,
   we have that, uniformly on compact subsets defined in
  \ef{z4}, as $k \to \iy$,
   \be
   \label{z6}
   g(s_k+s) \to h(s) \whereA h_s=\AAA h, \quad h=0
   \,\,\,\mbox{at} \,\,\,\xi=0, \quad h|_{\xi=+\iy}=1.
    \ee

Consider this {\em limit} (at $s=+\iy$) {\em equation} obtained
from \ef{z3}:
 \be
  \label{z7}
   \tex{
  h_s= \AAA h \equiv  h_{\xi\xi}+ \frac 12\, h_\xi \inB \re_+ \times \re_+,
  \quad h(0,s)=0, \quad
  h(+\iy,s)=1.
  }
   \ee
It is a  linear parabolic PDE in the unbounded domain $\re_+$,
governed by the operator
 $\AAA$ admitting a standard symmetric representation in a weighted space.
 Namely, we have:

 \begin{proposition}
  \label{Pr.91}
  {\rm (i)} \ef{z7} is a gradient system in a weighted
  $L^2$-space, and

  {\rm (ii)} for bounded orbits,  the $\o$-limit set $\O_0$ of \ef{z7} consists of a unique
  stationary profile
  \be
 \label{z12}
  \tex{
  g_0(\xi)=1-{\mathrm e}^{- \xi/2},
  }
 \ee
 and $\O_0$ is uniformly stable in the Lyapunov sense in a
 weighted $L^2$-space.
\end{proposition}

\noi{\em Proof.}
 As a 2nd-order  equation, \ef{z7}
 is written in a symmetric form,
  \be
  \label{z71}
  {\mathrm e}^{\xi/2} h_s= ({\mathrm e}^{\xi/2} h_\xi)_\xi,
   \ee
   and hence admits multiplication by $h_s$ in $L^2$ that yields
   a monotone Lyapunov function:
    \be
    \label{z72}
     \tex{
    \frac 12 \, \frac {\mathrm d}{{\mathrm d}s} \int {\mathrm
    e}^{\xi/2}(h_\xi)^2 \, {\mathrm d} \xi=- \int {\mathrm e}^{\xi/2}(h_s)^2\, {\mathrm d} \xi \le 0.
    }
    \ee
Note that, in \ef{z71}, the derivatives  $h_\xi$ and $h_s$ have to
have an exponential decay
 at infinity in order the seminorms involved to make sense. It is
 essential that the limit profile \ef{z12} perfectly suits both.

Thus,  the problem \ef{z3} is a perturbed {\em gradient system},
that allows  to pass to the limit $s \to +\iy$ by using power
tools of gradient system theory; see e.g., Hale \cite{Hale}.

\ssk

 (ii) For a given bounded orbit $\{h(s)\}$, denote
$h(s)=g_0+w(s)$, so that $w(s)$ solves the same equation \ef{z71}.
 Multiplying by $w(s)$ in $L^2$ yields
  \be
  \label{z73}
   \tex{
    \frac 12 \, \frac{{\mathrm d}}{{\mathrm d}s} \, \int{\mathrm
    e}^{\xi /2} w^2 \, {\mathrm d} \xi=- \int{\mathrm
    e}^{ \xi /2}(w_\xi)^2 \, {\mathrm d} \xi <0
    }
    \ee
for any nontrivial solutions, whence the  uniform stability
(contractivity) property. \quad $\qed$

 \ssk

Finally, we state the main stabilization result in the boundary
layer, which establishes the actual class of generic solutions we
are dealing with.

\begin{theorem}
 \label{Th.Gen1}
 {\rm (i)} There exists a class of solutions of the perturbed equation \ef{z3}, for
 which, in a weighted $L^2$-space and uniformly on compact
 subsets,
  \be
  \label{z10}
  g(\xi,s) \to g_0(\xi) \asA s \to +\iy.
  \ee

  {\rm (ii)} \ef{z10} is particularly  true for all nontrivial
  nonnegative 
   solutions.
 \end{theorem}

 \noi {\em Proof.} (i) Under given hypotheses, the uniform
 stability result in (ii) of Proposition \ref{Pr.91} implies  \cite[Ch.~1]{AMGV} that
 the $\o$-limit set of the asymptotically perturbed equation
 \ef{z3} is contained in that for the limit one \ef{z7}, which
 consists of the unique profile \ef{z12}.

 (ii) This follows from the construction, since then $\rho(s)$ in
 \ef{z2} can be chosen always positive. Then in the limit we are
 guaranteed to arrive at the gradient problem \ef{z6} admitting the unique
 uniformly stable stationary point \ef{z12}.
  \quad $\qed$

\section{Inner Region expansion: towards an ODE regularity
criterion}
 \label{S3.4}

\subsection{The Cauchy problem setting, eigenfunction expansion, and matching}

In Inner Region, we deal with the original rescaled problem
 \ef{ph4}.
Without loss of generality, again for simplicity of final, rather
technical and involved   calculations, we
 consider even solutions defined  for $y>0$ by
 assuming  the symmetry condition at the origin
   \be
   \label{ss111}
  v_y=0 \quad \mbox{at\,\,\, $y=0$}.
 \ee
As customary in classic PDE and
  potential theory (see e.g., Vladimirov \cite[\S~6]{Vlad71}),
  we extend
 $v(y,\t)$ by 0 beyond the boundary points, i.e., for $y > \var(\t)$:
  \be
  \label{a1}
 \hat v(y,\t) = v(y,\t)H(\var(\t)-y)=
  \left\{
   \begin{matrix}
   v(y,\t) \forA 0 \le y < \var(\t), \\
\,\, 0  \,\,\forA \,\,  y \ge  \var(\t),
\end{matrix}
 \right.
 \ee
 where $H$ is the Heaviside function.
Since $v=0$ on the lateral boundary $\{y= \var(\t)\}$, one can
check that, in the sense of distributions,
 \be
 \label{a2}
 \hat v_\t= v_\t H, \quad \hat v_y= v_y H, \andA
 \hat v_{yy}= v_{yy} H - v_{y}\big|_{y=\var}\d(y-\var).
  \ee
Therefore, $\hat v$ satisfies the Cauchy problem:
 \be
 \label{a3}
  \hat v_\t= \BB^* \hat v + v_{y}\big|_{y=\var(\t)}\d(y-\var(\t)) + \kappa(\hat v) \hat v
   \inB \re \times \re_+.
  \ee

Since, by construction, the extended solution  \ef{a1} is
uniformly bounded in $L^2_{\rho^*}(\re)$, we can use the
converging in the mean (and uniformly on compact subsets in $y$)
the eigenfunction expansion via the standard Hermite polynomials
given in \ef{psi**1} for $m=1$:
 \be
 \label{a4}
  \tex{
  \hat v(y,\t)= \sum_{(k \ge 0)} a_k(\t) \psi_k^*(y).
   }
   \ee
 Actually, as follows from BL-theory from Section \ref{S3.3},
 Theorem \ref{Th.Gen1}, that  the only possible solutions admitting matching with \ef{z10}
 possess  a
 constant in $y$ behaviour on compact subsets in $y$, i.e.,
 \be
 \label{ass12}
 \hat v(y,\t)= a_0(\t) \cdot 1 (1+o(1)) \asA \t \to + \iy.
  \ee
  Indeed, this ``$1$" well corresponds to the first Hermite
  polynomial $\psi_0^*(y) \equiv 1$ in \ef{a4}. Since $\l_0=0$ for
  this ``polynomial", the behaviour \ef{ass12} can be referred as
  to a ``center subspace" one for the operator $\BB^*$ in
  \ef{BB*N}, though we do not  use this fact at all.


Thus, by BL-theory establishing the boundary behaviour \ef{z2} for
$\t \gg 1$, which we state again: in the rescaled sense, on the
given compact subsets,
 \be
 \label{9}
  \tex{
 \hat v(y,\t) = \rho(s) g_0\big( \var^{2}(\t)(1- \frac
  y{\var(\t)}) \big)(1+o(1)).
  }
  \ee
 Overall, in the class of generic
 solutions satisfying the BL-expansion,
   we   concentrate on the  first Fourier  pattern associated with
 \be
 \label{a6}
 k=0: \quad \l_0=0 \andA \psi_0^*(y) \equiv 1 \quad \big(\psi_0(y) \equiv F(y), \,\,\,
  \mbox{the Gaussian
 (\ref{Fm1})}\big).
  \ee
 The corresponding normalization condition  is key for further
 projections:
  \be
  \label{a6N}
   \tex{
   \langle \psi_0, \psi^*_0 \rangle \equiv \int F=1.
   }
    \ee


\begin{proposition}
\label{Pr.Stable} Under the given assumptions: {\rm (i)} for
 solutions  in Theorem $\ref{Th.Gen1}${\rm (i)}, \ef{ass12} holds with
 $a_0(\t) > 0$ for $\t \gg 1$,
  and then the matching with the
 boundary layer behaviour in \ef{z2} requires
 \be
  \label{10}
   \tex{
   \frac {a_0(\t)}{\rho(s)} \to 1 \asA \t \to +\iy
   \LongA \rho(s) = a_0(\t)(1+o(1)).
   }
   \ee

{\rm (ii)} In particular, these are true for
 nontrivial
nonnegative 
 solutions.
 \end{proposition}

 \noi {\em Proof.} (i) follows from the construction of the
 boundary layer.

 (ii) This follows from Theorem \ref{Th.Gen1}(ii).
 \quad $\qed$

\ssk

 Thus,  projecting the PDE \ef{a3} onto the centre subspace of $\BB^*$ (i.e.,
 by multiplying in $L^2$ by $\psi_0(y)=F(y)$)
yields, for the leading mode $a_0(\t)$, the following ``ODE":
\be
   \label{a8}
   \tex{
   a_0'=  v_{y}\big|_{y=\var(\t)}  \psi_0(\var(\t))  +
\langle \kappa(\hat v)\hat v, \psi_0 \rangle.
    }
    \ee
 The convergence \ef{9}, which by a standard parabolic regularity is also true for the
 spatial derivatives, yields,
  as
 $\t \to +\iy$,
  \be
  \label{11}
  \tex{
 v_{y}\big|_{y=\var(\t)} = \rho(s) \var(\t) \g_1 (1+o(1)) = a_0(\t)
 \var(\t) \g_1(1+o(1)), \quad \g_1= g_0'(0)= \frac 12.
 }
   \ee

Finally, we need to estimate the last term in \ef{a8}: by
\ef{ass12}, using that $\kappa(a_0(\t)) \not =0$ for any $a_0(\t)
\ne 0$ via \ef{kap1}, there holds
 \be
  \label{11N}
  \tex{
  \langle \kappa(\hat v) \hat v, \psi_0 \rangle =  \langle \kappa(a_0) a_0, F
  \rangle (1+o(1))= \kappa(a_0) a_0(1+ o(1)).
  }
  \ee
Indeed, since $\int F=1$ for the Gaussian \ef{Fm1}, in the last
estimate, we have
 $$
 \tex{
  \int\limits_0^\var F(y)\, {\mathrm d}y \equiv \frac 12 -
  \int\limits_\var^\iy F(y) \, {\mathrm d}y= \frac 12- O\big( \frac 1
  \var \, {\mathrm e}^{-\var^2/4}\big) \asA \var=\var(\t) \to +\iy.
  }
  $$

   Thus, bearing in mind all above assumptions and estimates for generic patterns including
 \ef{ass12},  \ef{10}, \ef{9}, and \ef{11N},  we obtain
 the following asymptotic ODE for the first
expansion coefficient $a_0(\t) \ne 0$: as $ \t \to +\iy$,
 \be
 \label{12}
 \fbox{$
  \tex{
   \frac {a_0'(\t)}{a_0(\t)}=
    - \frac 1{4 \sqrt \pi}\,
   \var(\t)\, {\mathrm e}^{-\var^2(\t)/4}(1+o(1)) + \kappa(a_0(\t))(1+o(1)) \,
   .
   }
   $}
   \ee
One can see that, by assumptions \ef{kap1}, all the solutions of
the non-autonomous ODE \ef{12} are well defined for $\t \in [0,
+\iy)$. Moreover, by classic comparison/monotonicity results for
ODEs (S.A.~Chaplygin's theorem of 1920s \cite{Chap}), it follows
that, under the above hypotheses, solutions of \ef{12} satisfy:
 \be
 \label{comp1}
 a_0(0) \ne 0 \LongA a_0(\t) \ne 0 \quad \mbox{for all} \quad \t >
 0.
  \ee
 Therefore, we can always consider positive orbits:
  \be
  \label{comp2}
  a_0(\t)>0 \quad \mbox{for all} \quad \t \ge 0.
   \ee
  This makes our further asymptotic analysis easier.
 In particular, in view of \ef{comp1} and \ef{kap1}, we  can
 always
 omit all higher-order terms appeared via the above asymptotics.

\subsection{ODE regularity criterion}

It follows from \ef{12}, that a natural way to formulate a
regularity criterion for the parabolic PDE \ef{ir12} is to use the
``ODE language"\footnote{In fact, this is quite natural and
unavoidable: for semilinear PDEs, characteristic point regularity
depends on asymptotic properties of ODEs, i.e., regularity issues
for  infinite-dimensional dynamical systems are characterized by
1D ones. This reveals a sufficient and successful reduction of
dimensions.}:


\begin{theorem}
 \label{Th.kap1}
 $(${\bf ODE regularity criterion}$)$
In the parabolic problem \ef{ir12}, the origin $(0,0)$ is regular,
iff \, $0$\, is globally asymptotically stable for the ODE
\ef{12}, i.e., any solution of \ef{12} is global and satisfies
 \be
 \label{12n}
  \tex{
  a_0(\t) \to 0 \asA \t \to +\iy, \quad \mbox{i.e.,} \quad  \ln|a_0(\t)| \to -
  \iy.
  }
  \ee
  \end{theorem}

\noi{\em Proof.} (i) {\em Necessity.} Given any classic solution
$u(x,t)$ \ef{ir12}, one can always construct positive and negative
{\em barrier} solutions $u_\pm(x,t)$ such that
 \be
 \label{res1}
 u_-(x,t) \le u(x,t) \le u_+(x,t)  \inB Q_0
 \ee
 by standard comparison (Maximum Principle) arguments, \cite{Fr0}. Since, by Theorem \ref{Th.Gen1}(ii)
  and Proposition \ref{Pr.Stable}, such non sign-changing
 solutions $u_\pm(x,t)$ do obey our matched asymptotics, their
 positive (resp., negative)
 first Fourier coefficients satisfy the asymptotic ODE \ef{12} for
 $\t \gg 1$. Hence, by the BL-construction, \ef{12n} implies that
 $u_\pm(x,t) \to 0$ as $t \to 0^-$ uniformly, so, by comparison
 \ef{res1}, the same does an arbitrary $u(x,t)$.

 \ssk

 (ii) {\em Sufficiency by contradiction.} Let there exist a solution $\{\bar
 a_0(\t)\}$ of \ef{12n} (by \ef{comp1}, we may assume it to be
 positive) such that
  \be
  \label{res2}
 \tex{
   \limsup_{\t \to +\iy} \bar a_0(\t) >0.
   }
 \ee
 Then, by the ODE comparison, the same is true for solutions of \ef{12n} with
 arbitrarily large
 Cauchy data at $\t=0$, i.e., for any
  \be
  \label{res3}
  a_0(0) > \bar a_0(0).
   \ee
Therefore, there exists a sufficiently large positive solution
$u_+(x,t)$  of \ef{ir12}, whose first Fourier coefficient
satisfies \ef{12} and \ef{res3}, so the regularity  is violated by
\ef{res2}. \quad $\qed$

\ssk

For the heat equation \ef{heat1}, with $\kappa =0$, integrating
\ef{12} immediately yields
  \be
   \label{993}
    \fbox{$
     \kappa=0: \quad
\mbox{$(0,0)$ is regular iff} \quad
 \int\limits^\iy \var(\t)\, {\mathrm e}^{-
 \frac{\var^2(\t)}4}\, {\mathrm d}\t = +\iy,
  $}
  \ee
  which is indeed another equivalent form of Petrovskii's
  criterion \ef{RR1} (in Khinchin's form).




\subsection{Applications: further regularity results}
 \label{S4.6}

We now present a few corollaries of Theorem \ref{Th.kap1}, with
simpler and more traditional conditions of
regularity/irregularity.

First of all, it follows from the ODE \ef{12} (and actually is
true by comparison) that {\em negative} coefficients $\kappa(v)$
can ``improve" the regularity of $(0,0)$. Moreover, in this
simpler case, we find a condition, under which {\em any backward
parabola has a regular vertex}.

\begin{proposition}
\label{Pr.99}
 Let $\kappa(u)$ satisfy \ef{kap1} and let
  \be
  \label{kap222}
  \kappa(u)<0 \quad \mbox{for} \,\,\, u>0.
  \ee
Then, for any backward parabola $\pa Q_0$ with arbitrary $\var$'s
in \ef{ph1}, \ef{vv1}, the vertex $(0,0)$ is regular.
\end{proposition}

\noi{\em Proof.} It follows from \ef{12} that, for $\t \gg 1$,
 \be
 \label{kap44}
  \tex{
   \frac{a_0'}{a_0} \le - | \kappa(a_0)|(1+o(1)) \le - \frac 12 \,  | \kappa(a_0)|.
    }
    \ee
 Then, on integration, assuming, without loss of generality, that
$a_0(0)=1$,
 and checking an Osgood--Dini-type condition
 \be
 \label{kap33}
  \tex{
  \int\limits_{0^+} \frac{{\mathrm d}z}{z |\kappa(z)|}  = \iy,
   }
   \ee
 which obviously holds for the coefficients \ef{kap1},
  we have
 \be
 \label{kap55}
  \tex{
   \int\limits_{a_0(\t)}^1 \frac{{\mathrm d} z}{z|\kappa(z)|} \ge  \frac
   \t2
   \to +\iy \asA \t \to \iy.
   }
   \ee
Hence, \ef{kap33} reinforces \ef{12n} to hold. \quad $\qed$

\ssk

 Secondly, for {\em positive} coefficients $\kappa$, regularity can be
 destroyed.
 We first state the result establishing the conditions on
 monotone $\kappa(v)>0$, under which the nonlinear term changes
 regularity for the pure heat equation into the irregularity.


\begin{proposition}
 \label{Pr.reg2}
 Let $\kappa(u)$ satisfy \ef{kap1} and let
  \be
  \label{kap2}
  \kappa(u)>0 \quad \mbox{be increasing for} \,\,\, u>0.
  \ee
  Let
   \ef{993} be valid, i.e.,
  $(0,0)$ is regular for the heat
 equation \ef{heat1} for $N=1$.
 Denote by $\hat a_0(\t) \to 0$ as $\t \to \iy$
 the corresponding Fourier
 coefficient satisfying \ef{12} for $\kappa=0$:
  \be
  \label{a11}
    \hat a_0(\t)=\hat a_0(0) \, \eee^{- \frac 1{4 \sqrt \pi}\, \int_0^\t \var(s) \, {\mathrm e}^{-
   {\var^2(s)}/4} \, {\mathrm d}s} \forA \t \gg 1.
     \ee
Then the linear regularity criterion \ef{993} fails for the
semilinear problem \ef{ir12} and $(0,0)$ becomes irregular
provided that the nonlinearity $\kappa$ is such that
 \be
 \label{a12}
 \tex{
 \int\limits^{+\iy}\big[- \frac 1{4 \sqrt \pi}\, \var(\t) \eee^{-
 \frac{\var^2(\t)}4}+  \kappa(\hat a_0(\t))\big]\,
 {\mathrm d}\t > - \iy.
 }
 \ee
    \end{proposition}

\noi{\em Proof.}
 One can see that, in the present proof of a sharp estimate, one can
 omit both  $o(1)$-terms in \ef{12}, meaning that one can replace
 those by $1+\e$ and $1-\e$ with an $\e>0$ respectively and pass
 to the limit $\e \to 0^+$.

 As the first iteration of the full ODE \ef{12},
we have, for $\t \gg 1$,
 \be
  \label{b1}
   \tex{
    \frac {a_0'}{a_0} \ge - \frac 1{4 \sqrt \pi} \, \var(\t) \,
    \eee^{- \frac{\var^2(\t)}4} \LongA a_0(\t) \ge \hat a_0(\t).
    }
    \ee
  In view of \ef{kap2}, \ef{a12}, we then obtain via the second iteration of
  \ef{12}:
\be
 \label{b2}
  \tex{
  \frac{a_0'}{a_0} \ge - \frac 1{4 \sqrt{\pi}}\,\var(\t)
  \eee^{-\frac{\var^2(\t)}4}  +   \kappa(\hat a_0(\t)) \forA \t \gg 1.
  }
  \ee
 Integrating this yields, by \ef{12}, that $(0,0)$ is no more
 regular for such nonlinear coefficients $\kappa(v)$.
\quad $\qed$

\ssk

\begin{corollary}
 \label{Cor.1}
 Under  the conditions of Proposition $\ref{Pr.reg2}$,
 Petrovskii's backward parabola with the $2 \sqrt{\log
 \log}$-factor \ef{ph2}
 is no more a regular vertex of $Q_0$ for the semilinear problem \ef{ir12} provided
 that
  \be
  \label{a14}
   \tex{
   \kappa(v) \gg |\ln |v||^{\frac 13}\,\, {\mathrm e}^{-(3 \sqrt \pi |\ln |v||)^{2/3}}
    \asA v \to 0.
    }
 \ee
  \end{corollary}

Thus, \ef{a14} is the estimate, where the function in \ef{kap99}
comes from.

\ssk

\noi{\em Proof.} It follows from \ef{12} with $\kappa=0$ that the
function \ef{a11} reads
 \be
 \label{a15}
  \tex{
 \hat a_0(\t) \sim  {\mathrm e}^{-\frac 1{3 \sqrt \pi}\,(\ln
 \t)^{3/2}} \asA \t \to \iy.
 }
 \ee
Substituting \ef{a15} into \ef{a12} and changing the
variable $\hat a_0(\t)=v$ yields \ef{a14}. \quad $\qed$






\ssk

 Further iterating inequalities such as \ef{b2},
 one can
obtain other sufficient conditions of the origin irregularity.
 For instance, if the integral in \ef{a12} still diverges to
 $-\iy$,
 integrating \ef{b2} gives the next iteration estimate: for $\t
 \gg 1$
\be
 \label{kap51}
  \tex{a_0(\t) \ge \hat a_0^{(1)}(\t) \equiv a_0(0)\, \eee^{
  \int\limits_0^\t \big[- \frac 1{4 \sqrt \pi}\, \var(\eta)
{\mathrm e}^{-\var^2(\eta)/4} + \kappa \big(C_1 \,{\mathrm e}^{-
\frac 1{4 \sqrt \pi}\, \int_0^\eta \var(s) {\mathrm
e}^{-\var^2(s)/4}\, {\mathrm d}s} \big)\big]\,{\mathrm d} \eta}, }
\ee
 where $C_1>0$ is some constant.
 Then, the next iteration  leads to an awkward looking  inequality:
   \be
   \label{kap52}
    \begin{matrix}
    \frac {a_0'}{a_0} \ge - \frac 1{4\sqrt \pi}\, \var(\t)
    \eee^{-\var^2(\t)/4} \qquad\qquad\qquad
 \ssk\ssk \\
   + \kappa \Big(a_0(0)\, \eee^{
  \int\limits_0^\t \big[- \frac 1{4 \sqrt \pi}\, \var(\eta)
{\mathrm e}^{-\var^2(\eta)/4} + \kappa \big(C_1 \,{\mathrm e}^{-
\frac 1{4 \sqrt \pi}\, \int_0^\eta \var(s) {\mathrm
e}^{-\var^2(s)/4}\, {\mathrm d}s} \big)\big]\,{\mathrm d}
\eta}\Big).
 \end{matrix}
 \ee
 Integrating it gives an estimate of $a_0(\t) \ge \hat
 a_0^{(2)}(\t)$ for $\t \gg 1$
 from below to be used also for the purpose of the irregularity,
 if $a_0^{(2)}(\t) \not \to 0$ as $\t \to +\iy$. If this fails, we
 then apply
the third iteration of the ODE \ef{12} again leading to a sharper
estimate from below for the regularity, etc.

 Since the number of such iterations can increase without bound (and
hence the same do the numbers of exponents and corresponding
integrals in the argument of $\kappa(\cdot)$ in \ef{kap52}, etc.),
it seems inevitable that a single and a simply finite integral
criterion of irregularity, similar to the Petrovskii one \ef{993},
cannot be derived for the nonlinear dynamical system \ef{12} in
the maximal generality. In other words, the ODE criterion of
Theorem \ref{Th.kap1} is a right way to regularity issues and is
even optimal.

\ssk

In more general cases of equations in \ef{ir12}, where, in our
notations,
 \be
 \label{kap6}
  \kappa=\kappa(x,t,u,u_x),
   \ee
the derivation of matched asymptotics remains the same. The only
difference is that, in accurate estimating of the integral in the
last nonlinear term in \ef{a8}, we should take into account that
$v_y \approx 0$ in the whole inner region due to the ``centre
subspace expansion" \ef{ass12}, so actually we integrate there
$\kappa(\cdot, 0)$. But this term must also include integrals over
the boundary layers close to $y=\pm \var(\t)$, where the solution
$v$ and its derivative $v_y$  is sharply given by \ef{9} with the
matching condition \ef{10}. We do not perform these general and,
at the same time, rather straightforward and not that principal
computations here, and restrict our attention to a particular
model:

 \subsection{Equations with a gradient-dependent nonlinearity}
  \label{S4.7}

 Let us very briefly consider equation \ef{quad1}. Then, the first
 rescaling \ef{ph3} gives the equation
  \be
  \label{ss1}
  v_\t= \BB^* v+ \kappa(v) v (v_y)^2.
   \ee
It is easy to check that the BL-analysis yields the same
asymptotics as in \ef{9}, with a similar proof. However, the
eventual derivation of the 1D dynamical system for the first
Fourier coefficient $a_0(\t)$ is now different: the nonlinear term
is much weaker, since $v_y  \approx 0$ on the centre subspace
patterns, except a $\big(\frac 1{\var(\t)}\big)$-neighbourhood of
the boundary point $y=\var(\t)$. Overall, the nonlinear
perturbation in \ef{12} is estimated as follows:
 \be
 \label{ss3}
 \begin{matrix}
 J(a_0)=\langle \kappa(a_0) a_0^2 \big[g_0'\big(\var^2 \big(1- \frac
 y\var\big)\big)\big]^2 \,(- \var )^2, \,\, \psi_0(y)
 \rangle \ssk\ssk \\
 = {\kappa(a_0) a_0^2}{\var^2} \,  \frac 1{2 \sqrt{\pi}}\,
   \int\limits_0^\var \big[g_0'\big(\var^2 \big(1- \frac
 y\var\big)\big)\big]^2
   {\mathrm e}^{- y^2/4}\, {\mathrm d}y,
   \end{matrix}
 \ee
 where  $\psi_0=F$ given by \ef{Fm1}.
Using the BL-profile \ef{z12} and setting $z= \frac y \var$ yields
 \be
 \label{ss4}
  \begin{matrix}
 J(a_0)=  {\kappa(a_0)a_0^2 \var^3} \frac 1{8 \sqrt \pi} \,
  \int\limits_0^1 {\mathrm e}^{-\var^2(1-z)} {\mathrm e}^{-\var^2z^2/4}
  \, {\mathrm d} z
    \ssk\\
  = {\kappa(a_0)a_0^2 \var^3} \frac 1{8 \sqrt \pi} \,
   {\mathrm e}^{-\var^2} \int\limits_0^1 {\mathrm e}^{\var^2
   z(1-\frac z4)} \, {\mathrm d}z.
  \end{matrix}
  \ee
 Estimating roughly the last integral as follows:
  $$
   \tex{
\int\limits_0^1 {\mathrm e}^{\var^2
   z(1-\frac z4)} \, {\mathrm d}z \le {\mathrm e}^{ \frac3 4 \var^2},
    }
    $$
 we obtain the following approximate dynamical system for $a_0(\t)>0$:
 \be
 \label{ss5}
  \tex{
   \frac{a_0'}{a_0}  \le  - \frac 1{4 \sqrt \pi} \, \var(\t) \,{\mathrm
   e}^{- \var^2(\t)/4}+ \frac 1{8 \sqrt \pi} \, {\kappa(a_0)a_0^2 \var^3(\t)}
   {\mathrm e}^{-\var^2(\t)/4}+... \, .
    }
    \ee
    This is enough for us to prove that
  the nonlinear perturbation is now much weaker than that in \ef{12}:

\begin{proposition}
 \label{Pr.Weak}
  For \ef{ss5},
 Petrovskii's double log-function \ef{ph2} forms a regular vertex
 $(0,0)$ for any function $\kappa(u)$ satisfying \ef{kap1}.
  \end{proposition}

  \noi {\em Proof.} Assuming that the linear term is dominant that
  creates the behaviour \ef{a15}, one can check that, on this
  $\hat a_0(\t)$, the nonlinear term in \ef{ss5} is always
  negligible, so \ef{b2} follows. \quad $\qed$

  \ssk

\subsection{Backward paraboloid in $\ren$}

More carefully,  aspects of checking regularity of the vertex of a
{backward paraboloid} in $\ren$ was done in \cite{GalMazNSE},
where the authors applied  matching techniques to the
Navier--Stokes equations in $\re^3$.  Now we present  a few
comments.

For the $N$-dimensional case \ef{sem1N}, the lateral boundary of
the domain $Q_0$ in $\re^{N+1}$ is given by a {\em backward
paraboloid} of the form
 \be
 \label{parab1}
  \tex{
  \pa Q_0: \quad
 \sqrt{ \sum_{i=1}^N a_i |x_i|^2}= \sqrt{-t} \,\, \var(\t), \quad \t= -
  \ln (-t), \quad a_i>0, \,\,\, \sum a_i^2=1.
   }
   \ee
Then a boundary layer close to the rescaled (via \ef{z1}) boundary
 \be
 \label{vvv1}
  \tex{
\pa \hat Q_0: \quad  \sum a_i |z_i|^2=1,
  }
  \ee
  leads to a linear elliptic problem, which can be solved.
  Moreover, in the direction of the unit inward normal ${\bf n}$
  to $\pa \hat Q_0$, the boundary layer profile $g_0(\xi)$ remains
  one-dimensional depending on the single variable
 \be
 \label{norm1}
  \eta= \xi \cdot {\bf n},
  \ee
  so that $g_0=g_0(\eta)$ is still given by  \ef{z12}.
Therefore, in the expanding domain with the boundary
 \be
 \label{norm2}
  \tex{
 \pa \tilde Q_0(\t): \quad \sum a_i |y_i|^2= \var(\t) \to +\iy
 \asA \t \to +\iy,
  }
  \ee
  the BL-profile is expressed in terms of the distance function:
   \be
   \label{norm3}
   g_0(y,\t)=1- \eee^{-\frac 12 \var(\t) \, {\rm dist} \,\{y, \pa
   \tilde Q_0(\t)\}}.
    \ee
  This allows us to
  apply the same blow-up scaling and matching techniques.

  The final ODE for $a_0(\t)$ takes a similar to \ef{12} form, with
  $\var$ in the first term replaced by $\var^N$, in a full accordance to
  \ef{pp1}. However, the computations get more involved and
  further coefficients of this asymptotic ODE will essentially
  depend on the geometric shape of the backward paraboloid
  \ef{parab1} in a neighbourhood of its characteristic vertex
  $(0,0)$. However, final regularity conclusions remain
  approximately the same as for $N=1$, including both cases of
  nonlinearities in \ef{sem1N}.





\begin{appendix}
\section*{Appendix A. Hermitian spectral theory for operator pair
$\{\BB, \,\BB^*\}$}
 \label{S2}
 \setcounter{section}{1}
\setcounter{equation}{0}

\begin{small}


For the maximal generality and further applications,
 we describe the necessary
 spectral properties of the linear $2m$th-order differential
 operator in $\ren$
\begin{equation}
 \label{B1*}
  \tex{
 \BB^* = (-1)^{m+1} \D^{m}_y - \frac {1}{2m}
 \, y \cdot \n_y,
 }
 \end{equation}
and of its adjoint $\BB$ in the standard $L^2$-metric given by
  \begin{equation}
 \label{B1}
  \tex{
 \BB = (-1)^{m+1} \D_y^{m} + \frac {1}{2m}
 \, y \cdot \n_y +  \frac{N}{2m}\, I \quad (I \,\,\, \mbox{denotes the identity}).
 }
 \end{equation}
 Both operators occur after global and blow-up scaling respectively of solutions
 of the {\em poly-harmonic equation}
 \be
  \label{Lineq}
  u_t = - (-\D)^m u \inA.
  \ee
 Of course, for $m=1$, \ef{B1*} and \ef{B1}  are classic Hermite self-adjoint  operators
 with completely known spectral properties; see e.g., Birman--Solomjak \cite[pp.~44-48]{BS}.
 However, for any $m \ge 2$, both operators \ef{B1*} and \ef{B1},
 though looking very similar to those for $m=1$,
 {\em are not symmetric}  and do not admit a self-adjoint
extension, so we follow more recent paper \cite{Eg4} in presenting
necessary spectral results.
In what follows, we mainly must concentrate on the less known case
$m \ge 2$, naturally assuming that, for the classic self-adjoint
case $m=1$, we can borrow any result from several textbooks and/or
monographs.

\subsection{Fundamental solution, rescaled kernel, and first estimates}
 \label{Sect3}

 We begin with the necessary
  fundamental solution $b(x,t)$  of the linear poly-harmonic parabolic equation
\ef{Lineq},
 which is of standard similarity form and satisfies, in the
 sense of bounded measures:
  \be
  \label{1.3R}
   b(x,t) = t^{-\frac N{2m}}F(y), \quad y= x/t^{\frac 1{2m}} \quad
   \mbox{such that}
   \quad b(x,0^+)= \d(x),
 \ee
 where $\d(x)$ is Dirac's delta.
 The rescaled kernel $F=F(|y|)$ is then the unique radially symmetric solution of the elliptic
 equation with the operator \ef{B1}, i.e.,
  \begin{equation}
\label{ODEf}
 {\bf B} F \equiv -(-\Delta )^m F + \textstyle{\frac 1{2m}}\, y \cdot
\nabla F + \textstyle{\frac N{2m}} \,F = 0
 \,\,\,\, {\rm in} \,\, \ren,  \quad \mbox{with} \,\,\, \textstyle{\int F =
 1.}
\end{equation}

For $m=1$, $F$ is the classic positive Gaussian
 \be
 \label{Fm1N}
  \tex{
  F(y)= \frac 1{(4 \pi)^{N/2}} \, {\mathrm e}^{-|y|^2/4}\, >0 \inB
  \ren.
  }
  \ee
 For any $m \ge 2$, the rescaled kernel function $F(|y|)$ is   oscillatory as $|y| \to \infty$ and
satisfies the estimate (for $m=1$, this is trivial, with $\a=2$)
 \cite{EidSys, Fedor}
\begin{equation}
\label{es11} 
 |F(y)| < D\,\,  {\mathrm e}^{-d_0|y|^{\alpha}}
\,\,\,{\rm in} \,\,\, \ren, \quad \mbox{where} \,\,\,
\a=\textstyle{ \frac {2m}{2m-1}} \in (1,2),
\end{equation}
for some positive constants $D$ and $d_0$ depending on $m$ and
$N$.

\subsection{Sharp estimates in one dimension}

For further use in our regularity study, we need some sharp
estimates of the rescaled kernel, which we present for $N=1$,
where the regularity analysis gets also rather involved.
 Taking the Fourier
transform in \ef{ODEf} leads to the expression
 \be
 \label{FundSol}
  \tex{
 F(y) = \a_0  \int\limits_0^\infty {\mathrm e}^{-s^{2m}} \cos(s y)\, {\mathrm
 d}s,
 }
 \ee
 where $\a_0>0$ is the normalization constant,
 and, more precisely \cite{Eid},
\begin{equation}
\label{Eidf}
 \tex{
  F(y) = \frac 1 {\sqrt{2\pi}}\,  \int\limits_0^\infty {\mathrm
  e}^{-s^{2m}}
\sqrt{s|y|}\,\, J_{- \frac 1 2}(s|y|)\, {\mathrm d}s \quad
\mbox{in} \,\,\, \re, }
\end{equation}
where $J_\nu$ denotes Bessel's function.
   The rescaled kernel $F(y)$
satisfies \ef{es11}, where $d_0$ admits an explicit expression;
see below.
 Such optimal exponential estimates of the fundamental solutions
 of higher-order parabolic equations are well-known and were first
 obtained by Evgrafov--Postnikov (1970) and Tintarev (1982); see
 Barbatis \cite{Barb, Barb04} for key references and results.

As a crucial  issue for the further boundary point regularity
study, we will need a sharper, than given by \ef{es11}, asymptotic
behaviour of the rescaled kernel $F(y)$ as $y \to +\iy$. To get
that, we re-write the equation \ef{ODEf} on integration once as
 \be
 \label{i1}
 \tex{
 (-1)^{m+1}F^{(2m-1)} + \frac 1{2m} \, y F=0 \inB \re.
 }
 \ee
 Using standard classic WKBJ asymptotics, we substitute into \ef{i1}
 the function
  \be
  \label{i2}
  F(y) \sim y^{-\d_0} \, {\mathrm e}^{a y^\a} \asA y \to + \iy,
   \ee
   exhibiting two scales.
This gives the algebraic equation for $a$,
 \be
 \label{i3}
  \tex{
 (-1)^m (\a a)^{2m-1}= \frac 1{2m}, \andA \d_0=  \frac{m-1}{2m-1}>0\,
 .
 }
  \ee

 Note that the slow algebraically decaying factor $y^{-\d_0}$ in
\ef{i2} is available for any $m \ge 2$. For $m=1$, this algebraic
factor is absent  for the exponential positive Gaussian profile
 \be
 \label{Fm1}
  \tex{
   F(y)= \frac 1{2 \sqrt \pi}\, {\mathrm e}^{-y^2/4} \quad (m=N=1).
  }
  \ee


By construction, one needs to get the root $a$ of \ef{i3} with the
maximal ${\rm Re}\, a<0$. This yields (see e.g., \cite{Barb,
Barb04} and \cite[p.~141]{GSVR})
 \be
 \label{i4}
  \tex{
 a= \frac{2m-1}{(2m)^\a} \big[-\sin\big( \frac{\pi}{2(2m-1)}\big) +
 \ii \cos\big( \frac{\pi}{2(2m-1)}\big)\big] \equiv -d_0 + \ii b_0
 \quad (d_0>0).
 }
 \ee
Finally, this gives the following double-scale asymptotic of the
kernel:
 \be
 \label{i5}
  \tex{
  F(y) =
   y^{-\d_0} \, {\mathrm e}^{-d_0 y^\a} \big[ C_1 \sin (b_0 y^\a)+
   C_2 \cos (b_0 y^\a)\big]+... \asA y \to + \iy ,
   }
   \ee
 where $C_{1,2}$ are real constants, $|C_1|+|C_2| \not = 0$.
 In \ef{i5}, we present the first two leading terms from the
 $m$-dimensional bundle of exponentially decaying asymptotics.

In particular, for the linear bi-harmonic operator in \ef{sem4}
($N=1$),
 we have
 \be
 \label{i6}
  \tex{
  m=2: \quad \a= \frac 43,  \quad d_0=3 \cdot 2^{-\frac{11}3},
 \quad b_0=3^{\frac 32} \cdot 2^{-\frac{11}3},
  \andA \d_0= \frac
  13.
   }
   \ee

\subsection{The discrete real spectrum and eigenfunctions of
 $\BB$}


Both linear operators ${\bf B}$ and  the corresponding adjoint
operator ${\bf B}^*$
 should be considered  in weighted $L^2$-spaces with the
weight functions induced by the exponential estimate of the
rescaled kernel (\ref{es11}). We again more  concentrate on the
non-self-adjoint case $m \ge 2$, and refer to \cite{BS} for the
classic one $m=1$.

 For $m \ge 2$,
we consider $\BB$ in the weighted space $L^2_\rho(\ren)$ with the
exponentially growing weight function
 \be
  \label{rho44}
  \rho(y) = {\mathrm e}^{a |y|^\a}>0 \quad {\rm in} \,\,\, \ren,
 \ee
  where $a \in (0,  2d_0)$ is a fixed
constant.
 We next
introduce a standard  Hilbert (a weighted Sobolev) space of
functions $H^{2m}_{\rho}(\ren)$ with the inner product and the
induced  norm
\[
 \tex{
 \langle v,w \rangle_{\rho} = \int\limits_{\ren} \rho(y) \sum\limits_{k=0}^{2m}
 D^{k}_y
 v(y) \, \overline {D^{k}_y w(y)} \,{\mathrm d} y, \quad
\|v\|^2_{\rho} = \int\limits_{\ren} \rho(y) \sum\limits_{k=0}^{2m}
|D^{k}_y
 v(y)|^2 \, {\mathrm d} y.
  }
\]
Then $H^{2m}_{\rho}(\ren) \subset L^2_{\rho}(\ren) \subset
L^2(\ren)$, and  $\BB$ is a bounded linear operator from $
H^{2m}_{\rho}(\ren)$ to $ L^2_{\rho}(\ren)$. The necessary
spectral properties of the operator $\BB$ are as follows
\cite{Eg4}:

\begin{lemma}
\label{lemspec}
 {\rm (i)}  The spectrum of $\BB$
comprises real simple eigenvalues only,
 \begin{equation}
\label{spec1}
 \tex{
 \sigma(\BB)=
\big\{\lambda_\b = -\frac k{2m}, \,\, k= |\b|= 0,1,2,...\big\}.
 }
\end{equation}
 {\rm (ii)} The eigenfunctions $\psi_\b(y)$ are given by
\begin{equation}
\label{eigen} \psi_\beta(y) =\textstyle{\frac{(-1)^{|\b|}}{\sqrt
{\b !}}} \, D^\beta F(y), \quad \mbox{for any} \,\,\,|\b|=k
\end{equation}
and form a complete set
 in $L^2({\re})$ and in $L^2_{\rho}({\re})$.

\noi {\rm (iii)} The resolvent $(\BB-\lambda I)^{-1}$
for  $\lambda \not \in \sigma(\BB)$ is a compact integral operator
in $L^2_{\rho}(\ren)$.
\end{lemma}


By Lemma \ref{lemspec}, the   centre and stable subspaces of $\BB$
are given by
$
E^c = {\rm Span}\{\psi_0= F\}$, $E^s = {\rm Span}\{\psi_\b, \,
|\b|>0\}$.
 Note also that the operators $\BB$
  has the zero Morse index, i.e., no
eigenvalues have positive real part.
 In the classic
Hermite case $m=1$ (the only self-adjoint case),
 the spectrum is again given by \ef{spec1} and the eigenfunction formula
\ef{eigen} with the rescaled kernel \ef{Fm1} generates standard
Hermite polynomials; see \cite[p.~48]{BS} for a full spectral
account for the operator $\BB$.

\subsection{The polynomial eigenfunctions of the operator $\BB^*$}

We now consider the adjoint operator (\ref{B1*}) in the weighted
space $L^2_{\rho^*}(\ren)$ ($\langle \cdot, \cdot
\rangle_{\rho^*}$ and $\|\cdot\|_{\rho^*}$ are the inner product
and the norm)
 with the   ``adjoint"  exponentially decaying weight
function
  \begin{equation}
\label{rho2}
 \tex{
 \rho^*(y) \equiv \frac  1 {\rho(y)} = {\mathrm e}^{-a|y|^{\a}} > 0.
 }
\end{equation}
 We ascribe to $\BB^*$ the domain
 $H^{2m}_{\rho^*}(\ren)$, which is dense in $L^2_{\rho^*}(\ren)$, and
 then
$$
 \BB^*: \,\, H^{2m}_{\rho^*}(\ren) \to L^2_{\rho^*}(\ren)
$$
 is a bounded linear operator.  $\BB$ is adjoint
 to $\BB^*$ in the usual sense: denoting by $\langle \cdot,\cdot \rangle $  the
inner product in the dual space $L^2(\ren)$, we have
  \begin{equation}
 \label{Badj1}
 \langle \BB v, w \rangle =  \langle v, \BB^* w \rangle
 \quad \mbox{for any} \,\,\, v \in H^{2m}_\rho(\ren) \andA
 w \in H^{2m}_{\rho^*}(\ren).
 \end{equation}
The eigenfunctions of $\BB^*$ take a particularly simple
polynomial form and are as follows:


 \begin{lemma}
\label{lemSpec2}
 {\rm (i)} $ \sigma(\BB^*)=\s(\BB)$.

 \noi{\rm (ii)} The eigenfunctions  $\psi^*_\b(y)$ of $\BB^*$ are polynomials
in $y$ of the degree $|\b|$ given by
 \begin{equation}
 \label{psi**1}
  \tex{
 \psi_\b^*(y) = \frac 1{\sqrt{\beta !}}
 \Big[ y^\b + \sum_{j=1}^{[|\b|/2m]} \frac 1{j !}(-\Delta)^{m j} y^\b \Big]
 }
 \end{equation}
and form a complete subset  in $L^2_{\rho^*}(\ren)$.

  \noi {\rm (iii)}
$\BB^*$ has a compact resolvent $(\BB^*-\lambda I)^{-1}$ in
$L^2_{\rho^*}(\ren)$ for $\lambda \not \in \sigma(\BB^*)$.
\end{lemma}



Of course, for $m=1$, \ef{psi**1} yields standard Hermite
polynomials, so, for $m \ge 2$, we call \ef{psi**1} {\em
generalized Hermite polynomials}. The bi-orthonormality condition
holds:
 \be
 \label{Ortog}
\langle \psi_\b, \psi_\g^* \rangle = \d_{\b\g}.
 \ee

\ssk

\noi{\bf Remark on closure.} This is an important issue for using
eigenfunction expansions of solutions.
 Firstly,
 as is well-known, for $m=1$, the sets of eigenfunctions are complete and closed in the
corresponding spaces, \cite{BS}.

Secondly, for $m \ge 2$, one needs some extra speculations.
Namely, using (\ref{Ortog}), we can introduce the subspaces of
eigenfunction expansions and begin with the operator $\BB$. We
denote by $\tilde L^2_\rho$ the subspace of eigenfunction
expansions $v= \sum c_\b \psi_\b$ with coefficients $c_\b =
\langle v, \psi^*_\b \rangle$ defined as the closure of the finite
sums $\{\sum_{|\b| \le M} c_\b \psi_\b\}$ in the norm of
$L^2_\rho$. Similarly, for the adjoint operator $\BB^*$, we define
the subspace $\tilde L^2_{\rho^*} \subseteq L^2_{\rho^*}$. Note
that since the operators are not self-adjoint and the
eigenfunction subsets are not orthonormal, in general, these
subspaces can be different from
$
 L^2_{\rho}$ and $L^2_{\rho^*}$, and
the equality is guaranteed in the self-adjoint case $m=1$,
$a=\frac 1 4$ only. For $m \ge 2$, in the above subspaces obtained
via suitable closure, we can apply standard eigenfunction
expansion techniques as in the classic self-adjoint case $m=1$.

\ssk

  For $m=2$ and $N=1$ (this simpler case will be treated in greater
  detail),
  the first ``adjoint" generalized Hermite polynomial eigenfunctions
  are:
 \be
 \label{psi44}
  \begin{matrix}
 \psi_0(y) = 1, \quad \psi_1(y)=y, \quad \psi_2(y) = \frac 1{\sqrt{2}}\,  y^2,
 \quad \psi_3(y) = \frac 1{\sqrt 6}\, y^3,\qquad \ssk\ssk\\
  \psi_4(y)
=  \frac 1 {\sqrt{24}}\, (y^4 + 24),\,
 \, \psi_5(y)= \frac 1{2 \sqrt{30}}\,(y^5+ 120\, y), \,\,
 \psi_6(y) =\frac 1 {12\sqrt{5}}\, (y^6 +360y^2), \quad
 \end{matrix}
 \ee
 etc.,
  with the corresponding eigenvalues $0$, $-\frac 14 $, $-\frac 12$, $-\frac 34$, $-1$, $- \frac 54$,
    $- \frac 32$, etc.


%


\end{small}
\end{appendix}

\begin{appendix}
\section*{Appendix B. Semilinear bi-harmonic equations}
 \label{S4}
 \setcounter{section}{2}
\setcounter{equation}{0}

\begin{small}



\subsection{Regularity problem setting}

 Here, we show how our approach can be extended to
higher-order PDEs, e.g., for the {\em semilinear bi-harmonic
equations} having similar nonlinearities, with also zero Dirichlet
conditions on $\partial Q_0$ and bounded initial data $u_0$ in
$Q_0 \cap \{t=-1\}$:
 \be
 \label{sem4}
 u_t= -u_{xxxx} +
\left\{
 \begin{matrix}
 \frac{\kappa(u) u}{(-t)}\ssk\\
 \kappa(u)u ( u_x)^4
 \end{matrix}
 \right.
  \inB Q_0,  \quad u=u_x=0 \,\,\, \mbox{on} \,\,\, \pa Q_0, \quad u(x,-1)=u_0(x).
  \ee
Then, after a proper similar matching with a boundary layer, we
again arrive a  nonlinear dynamical system viewed as a ``centre
subspace" approximation of solutions in the space of generalized
Hermite polynomials as eigenfunctions of a rescaled
non-self-adjoint operator. We also discuss the regularity problems
for {\em backward paraboloids} $\pa Q_0$ in $\ren \times [-1,0)$,
where the IBVP reads
  \be
  \label{sem4N}
 u_t= -\D^2 u +
 \left\{
 \begin{matrix}
 \frac{\kappa(u) u}{(-t)} \ssk\\
 \kappa(u)u |\n u|^4
 \end{matrix}
 \right.
 \inB Q_0,\quad u= \tex{\frac {\pa u}{\pa {\bf n}}=0} \,\,\, \mbox{on} \,\,\, \pa
 Q_0,
 \ee
 where ${\bf n}$ is the unit inward normal to the smooth boundary
 of the domain
 $Q_0 \cap\{t\}$.
 Further extensions to $2m$th-order parabolic PDEs are also
  discussed.

\ssk

 Thus, we now show that a similar sequence of mathematical
 transformations can be performed for the fourth-order semilinear bi-harmonic
 equations \ef{sem4}.

 \subsection{IBVP}

We again fix $N=1$, i.e., consider \ef{sem4} with  the simplest
nonlinearity \ef{3.1}, leading to the IBVP
 \be
 \label{ir124}
  \left\{
 \begin{matrix}
   u_t=-u_{xxxx} + \frac 1{(-t)}\, \kappa(u)u \quad \mbox{in}
 \quad Q_0=\{|x|<R(t), \,\,\, -1 < t<0\},
 \ssk\ssk\\
  u=u_x=0 \atA x= \pm R(t), \quad -1 \le t <0,\quad\,\,\,
  \,\,\,\,\,\,\,\qquad\qquad\,\qquad
    \ssk\ssk\\
  u(x,0)=u_0(x) \onA [-R(-1),R(-1)],\qquad\quad\,\,\,\,\,\,\,\,
  \qquad\qquad\,\qquad
   \end{matrix}
    \right.
  \ee
where $u_0(x)$ is   bounded and satisfies $u_0=u_0'=0$ at $x=\pm
R(-1)$.

\subsection{Slow growing factor $\var(\t)$}

Similar to \ef{PP1}, we  assume that
 \be
 \label{ph14}
 R(t)= (-t)^{\frac 14}\, \var(\t) \whereA \t= - \ln(-t) \to + \iy
 \asA t \to 0^-.
  \ee
 Here, the main scaling  factor  $(-t)^{1/4}$ naturally comes from the bi-harmonic
 kernel variables (see \ef{ph3} and \ef{1.3R}),
and $\var(\t)>0$ is again a slow growing function
 satisfying \ef{vv1}. For ``shrinking backward parabolae" with
  $$
 \var(\t),\,\,\, \var'(\t) \to 0 \asA \t \to \iy,
  $$
   the regularity in the linear case $\kappa=0$ was  proved
 by Miha${\rm\check{i}}$lov in 1963
 \cite{Mih63I, Mih63II}; in a certain sense, this extended the Gevrey-like result
 \ef{Gev1} for $m=1$; see \ef{Gev1Mih}.

\subsection{First kernel scaling}

 By \ef{ph1}, we
perform the similarity scaling
 \be
 \label{ph34}
  \tex{
  u(x,t)=  v(y,\t) \whereA y=\frac x{(-t)^{
  1/4}}.
  }
   \ee
   The rescaled function $v(y,\t)$ solves the
   rescaled IBVP
 \be
 \label{ph44}
  \left\{
  \begin{matrix}
  v_\t= \BB^* v +\kappa(v)v\equiv - v_{yyyy}- \frac 14 \, y v_y +\kappa(v)v
  \inB Q_0=\{|y| < \var(\t), \,\,\, \t>0\}, \ssk\ssk\\
  v=v_y=0 \atA y= \pm \var(\t), \,\,\t \ge 0,
  \qquad\qquad\qquad\qquad\qquad\,\qquad\qquad\,\qquad\quad
  \ssk\ssk\\
  v(0,y)=v_0(y) \equiv u_0(y) \onA
  [-R(-1),R(-1)].\qquad\quad\,\,\,\,\,\,\,\,\,\,\,\,\qquad\qquad\,\qquad\quad
   \end{matrix}
   \right.
  \ee

 \subsection{Boundary layer}
   \label{S3.34}

 Sufficiently close to the lateral boundary of $Q_0$, we
 naturally  introduce the variables
  \be
  \label{z14}
   \tex{
   z= \frac y{\var(\t)}, \,\, v(y,\t)=w(z,\t)
    \LongA w_\t= - \frac 1{\var^4} \, w_{zzzz} - \frac 14\, z w_z +
    \frac {\var'}\var \, z w_z + \kappa(w)w.
    }
    \ee
 The BL-variables  now read
  \be
  \label{z24}
  \tex{
  \xi= \var^{\frac 43}(\t)(1-z), \quad \var^{\frac 43}(\t) {\mathrm d} \t={\mathrm d}s,
   \andA w(z,\t)= \rho(s) g(\xi,s),
  }
  \ee
  where $\rho(s)$ is a slow varying  function, for which
  eventually \ef{10} will hold by matching.

Substituting into \ef{z14} yields the  perturbed equation
 \be
 \label{z34}
  \begin{matrix}
  g_s= \AAA g - \frac 14\,  \frac 1{\var^{4/3}} \, \xi g_\xi - \frac
  {\var'_\t}{\var}\, \big(1- \frac \xi {\var^{4/3}}\big) g_\xi \ssk\ssk\\
- \frac 43\,  \frac{\var'_\t}{\var^{1/3}} \,
   \xi g_\xi -  \frac {\rho'_s}{\rho} \,
  g + \frac 1{\var^{4/3}}\, \kappa(\rho g)\, g, \quad  \mbox{where} \quad  \AAA g= - g^{(4)} + \frac 14\, g'.
    \end{matrix}
   \ee
In this boundary layer,  we are looking for a generic
 pattern of the behaviour described by \ef{z34} on compact subsets
 near the lateral boundary,
  \be
  \label{z44}
  |\xi| = o\big(\var^{-\frac 43}(\t)\big)
   \LongA |z-1| = o\big(\var^{-\frac 83}(\t)\big) \asA \t \to +
   \iy.
   \ee

 We next pose the same asymptotic behaviour \ef{z5} at infinity.
 Assuming that, by  \ef{z24},
  the rescaled orbit $\{g(s),\,\, s>0\}$ is uniformly bounded, by
   parabolic theory \cite{EidSys}, we can again pass to the
  limit in \ef{z34} in the asymptotically small perturbations, along a subsequence $\{s_k\} \to +\iy$.
  Therefore,
   uniformly on compact subsets defined in
  \ef{z44}, as $k \to \iy$,
   \be
   \label{z64}
   g(s_k+s) \to h(s) \whereA h_s=\AAA h, \quad h=h_\xi=0
   \,\,\,\mbox{at} \,\,\,\xi=0, \quad h|_{\xi=+\iy}=1.
    \ee
The {limit equation} obtained from \ef{z34},
 \be
  \label{z74}
   \tex{
  h_s= \AAA h \equiv - h_{\xi\xi\xi\xi}+ \frac 14\, h_\xi
  }
   \ee
is again a standard linear parabolic PDE in the unbounded domain
$\re_+$, though now it is governed by a non self-adjoint operator
 $\AAA$.
 Actually, we need to show that, in an appropriate weighted
 $L^2$-space if necessary and under the hypothesis \ef{z5},
 the stabilization holds, i.e.,
  the $\o$-limit set of the orbit $\{h(s)\}_{s>0}$ consists of a
  single
 equilibrium: as $s \to +\iy$,
  \be
  \label{z104}
   \left\{
   \begin{matrix}
  h(\xi,s) \to g_0(\xi) \whereA \AAA g_0=0 \,\,\, \mbox{for}
  \,\,\,\xi>0, \ssk\ssk \\
   g_0=g_0'=0 \quad \mbox{at} \quad \xi=0, \quad
  g_0(+\iy)=1.\qquad\,\,
  \end{matrix}
  \right.
   \ee
The characteristic equation for the linear operator $\AAA$ yields
 \be
 \label{z114}
  \tex{
  -\l^4 + \frac 14\, \l=0 \LongA \l_1=0 \andA \l_{2,3}= \frac 1{4^{1/3}}
  \big(-\frac 12 \pm \ii \frac{\sqrt{3}}2\big).
  }
  \ee
  This gives the unique solution of \ef{z104},  shown in Figure
  \ref{F1},
 \be
 \label{z124}
  \tex{
  g_0(\xi)=1-{\mathrm e}^{- \frac \xi{2^{5/3}}} \,\big[ \cos \big( \frac
  {\sqrt{3} \, \xi}{2^{5/3}}\big) + \frac 1{\sqrt 3} \, \sin \big( \frac
  {\sqrt{3}\, \xi}{2^{5/3}}\big) \big].
  }
  \ee

\begin{figure}
\centering
\includegraphics[scale=0.65]{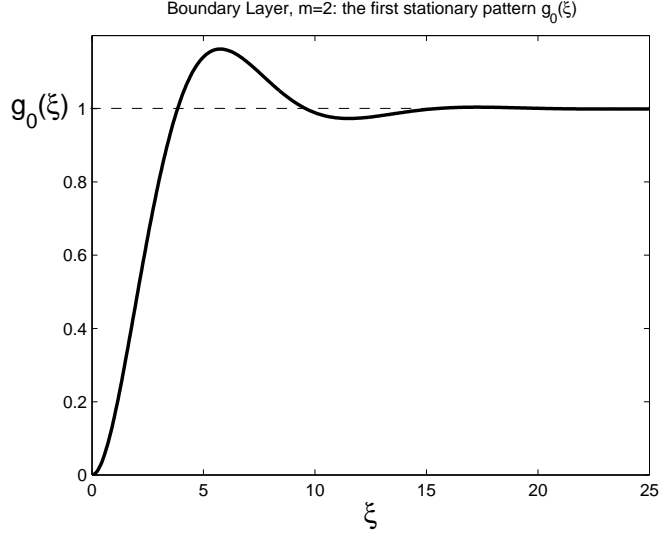} 
\vskip -.3cm \caption{\small \cite{GalPet2m} The unique stationary
solution $g_0(\xi)$ of the problem \ef{z104}.}
   \vskip -.3cm
 \label{F1}
\end{figure}

It turns out that the limit problem \ef{z74} possesses a number of
strong gradient and contractivity properties. Namely setting by
linearization
 \be
 \label{g65}
 \tex{
 h(s)=g_0 + w(s) \LongA w_s= \AAA w \equiv - w_{\xi\xi\xi\xi} +
 \frac 14 \, w_\xi, \quad w=w_\xi=0 \,\,\, \mbox{at} \,\,\, \xi=0,
 }
 \ee
we arrive at the following (cf. Proposition \ref{Pr.91} for
$m=1$):

\begin{proposition}
 \label{Pr.Grad4}
 {\rm (i)} \ef{g65} is a gradient system in $L^2$, and

 {\rm (ii)} in the given class of solutions, the $\o$-limit set
 $\O_0$ of \ef{g65} consists of the origin only and it is
 uniformly stable.
 \end{proposition}

 \noi{\em Proof.} (i) One can see that \ef{g65} admits a
 monotone Lyapunov function obtained by multiplying
 $w_{\xi\xi}$ in $L^2$:
  \be
  \label{g66}
 \tex{
   \frac 12\, \frac {\mathrm d}{{\mathrm d}s} \int (w_\xi)^2= -
   \int(w_{\xi\xi\xi})^2 \le 0.
   }
    \ee
Hence, (ii) also follows. \quad $\qed$

\ssk

Thus, quite similar to the second-order case, under given
assumptions, we can pass to the limit $s \to +\iy$ along any
sequence in the perturbed gradient system \ef{z34}. Then, again
similarly to $m=1$, the uniform stability of the stationary point
$g_0$ in the limit autonomous system \ef{z74} in a suitable metric
 guarantees that the asymptotically small perturbations do not
affect the omega-limit set; see \cite[Ch.~1]{AMGV}.
 However,
at this moment, we cannot avoid the following convention, which
for $m =2$ is much more key than for $m=1$, where the Maximum
Principle makes this part of the analysis simpler, at least, for
nonnegative or non-positive solutions (but for others of changing
sign, this remains   necessary). Actually, the convergence
\ef{z64} and \ef{z104} for the perturbed dynamical system \ef{z34}
should be considered as the main hypothesis,
{\em characterizing the class of generic patterns} under
consideration (and then \ef{z5} is its partial consequence).
 Since the positivity (negativity) is not an invariant property for
 bi-harmonic equations, a more clear characterization of this class of
 generic patterns is difficult. It seems that a correct language
 of doing this (in fact, for both cases $m=1$ and $m \ge 2$)
  is to reinforce  a ``centre subspace behaviour" as in
 \ef{ass12}, rather than other (possibly, ``stable") ones.
  Or, equivalently (and even more
 solidly mathematically), to impose the BL-behaviour  \ef{z104}.

Finally, we
 summarize these conclusions as follows:

\begin{proposition}
 \label{Pr.g04}


Under the given hypothesis and conditions, the
 problem $\ef{z34}$ admits a family of solutions $($called  generic$)$
  satisfying $\ef{z104}$.

 \end{proposition}

Such a definition of generic patterns looks rather
non-constructive, which is unavoidable  for higher-order parabolic
PDEs without positivity and order-preserving features. However, we
expect that \ef{z104} occurs for ``almost all" solutions.

\subsection{Inner region analysis: towards the dynamical system}
 \label{S3.44}

As usual, in the Inner Region, we treat the original rescaled
problem
 \ef{ph44}.
For simplicity of calculations, we again
 consider symmetric solutions defined  for $y>0$ by
 assuming  the symmetry at the origin:
   \be
   \label{ss1114}
  v_y=v_{yyy}=0 \quad \mbox{at\,\,\, $y=0$}.
 \ee
  We next extend
 $v(y,\t)$ by 0 for $y > \var(\t)$ and use the change \ef{a1}.
Since $v=v_y=0$ on the lateral boundary $\{y= \var(\t)\}$, one can
check that, in the sense of distributions,
 \be
 \label{a24}
  \begin{matrix}
 \hat v_\t= v_\t H, \quad \hat v_y= v_y H, \quad
 \hat v_{yy}= v_{yy} H, \ssk\ssk\\
 \hat v_{yyy}= v_{yyy} H - v_{yy}\big|_{y=\var}\d(y-\var),\ssk\ssk
 \\
 \hat v_{yyyy}= v_{yyyy} H - v_{yyy}\big|_{y=\var}\d(y-\var)-
 v_{yy}\big|_{y=\var}\d'(y-\var).
 \end{matrix}
  \ee
Therefore, $\hat v$ satisfies the following equation:
 \be
 \label{a34}
  \hat v_\t= \BB^* \hat v - v_{yyy}\big|_{y=\var}\d(y-\var)-
 v_{yy}\big|_{y=\var}\d'(y-\var) + \kappa(\hat v) \hat v \inB \re_+\times \re_+.
  \ee
Since such an  extended solution orbit \ef{a1} is uniformly
bounded in $L^2_{\rho^*}(\re)$, we  use  the eigenfunction
expansion via the generalized Hermite polynomials \ef{psi**1}:
 \be
 \label{a44}
  \tex{
  \hat v(y,\t)= \sum_{(k \ge 0)} a_k(\t) \psi_k^*(y).
   }
   \ee
   Substituting \ef{a44} into \ef{a34} and using the bi-orthonormality
   property \ef{Ortog} yields a dynamical system: for  $k=0,1,2,...\, ,$
   \be
   \label{a54}
   \tex{
   a_k'= \l_k a_k -  v_{yyy}\big|_{y=\var(\t)} \langle
   \d(y-\var(\t)), \psi_k \rangle - v_{yy}\big|_{y=\var(\t)}
   \langle \d'(y-\var), \psi_k \rangle + \langle \kappa(\hat v)\hat v, \psi_k \rangle,
    }
    \ee
    where $\l_k= -\frac k 4$ by \ef{spec1}. Here,
    $\l_k<0$ for all $k \ge 1$. More importantly, the
    corresponding eigenfunctions $\psi_k^*(y)$ are unbounded polynomials and are not
    monotone for $k \ge 1$ according to \ef{psi44}. Therefore,
    regardless proper asymptotics given by \ef{a54}, these inner
    patterns cannot be matched with the BL-behaviour such as
    \ef{z5}, and demand other matching theory. However, these are not
    generic, so we skip them.


   Thus, we  concentrate on the ``maximal" first Fourier generic pattern associated with
 \be
 \label{a64}
 k=0: \quad \l_0=0 \andA \psi_0^*(y) \equiv 1 \quad
 \big(\psi_0(y)=F(y)\big),
  \ee
  which corresponds to a
 ``centre subspace behaviour" \ef{ass12} for
 the equation \ef{a54}, which can be treated as
 another characterization of our class of generic patterns.
   The equation for $a_0(\t)$ is:
\be
   \label{a84}
   \tex{
   a_0'= -  v_{yyy}\big|_{y=\var(\t)}  \psi_0(\var(\t))  + v_{yy}\big|_{y=\var(\t)}
    \psi'_0(\var(\t)) + \langle \kappa(a_0)a_0, \psi_0 \rangle +..\,.
    }
    \ee

Next, we use the boundary behaviour \ef{z24}, \ef{z104} for $\t
\gg 1$, which for convenience we state again: in the rescaled
sense, on the given compact subsets,
 \be
 \label{94}
  \tex{
  v(y,\t) = \rho(s) g_0\big( \var^{\frac 43}(\t)(1- \frac
  y{\var(\t)}) \big)+...\,,
  }
  \ee
  where $g_0$ is as in \ef{z124}.
  Then,
 by the matching of both Regions for such
 generic patterns, \ef{10} must remain valid.
 Therefore, by  \ef{94}, which by a standard parabolic regularity is also true for the
 spatial derivatives, we have that, as $\t \to +\iy$,
  \be
  \label{114}
   \begin{matrix}
 v_{yy}\big|_{y=\var(\t)} \to \rho(s) \var^{\frac 23}(\t) \g_1 \to a_0(\t)
 \var^{\frac 23}(\t) \g_1 \whereA \g_1= g_0''(0)=2^{-\frac 43},\qquad \ssk\ssk\\
v_{yyy}\big|_{y=\var(\t)} \to -\rho(s) \var(\t) \g_2 \to - a_0(\t)
 \var(\t) \g_2 \whereA \g_2= g_0'''(0)= -\frac 14. \qquad
  \end{matrix}
   \ee

Finally, for such generic patterns,  we arrive at the asymptotic
ODE for the first Fourier coefficient:
 \be
 \label{124}
  \fbox{$
  \tex{
   \frac {a_0'}{a_0}=  G_4(\var(\t),\kappa) \equiv \g_2\var(\t) \psi_0(\var(\t))+ \g_1
   \var^{\frac 23}(\t) \psi_0'(\var(\t))+ \kappa (a_0)+...  \forA \t \gg
   1,
   }
   $}
   \ee
   where, as usual, we omit higher-order terms relative all those
   remaining. Note that, for this ODE, the properties \ef{comp1}
   and \ef{comp2} remain valid by comparison.

 \subsection{ODE regularity criterion and further applications}

In general, the criterion of regularity \ef{12n} remains the same.
However, it now reads:

\begin{theorem}
 \label{Th.kap14}
 $(${\bf ODE regularity criterion}$)$
In the fourth-order parabolic problem \ef{ir124}, the origin
$(0,0)$ is regular in the class of generic solutions, iff any
solution of the ODE \ef{124} satisfies \ef{12n}.
  \end{theorem}

Recall that by generic solutions we mean those that obey the
boundary layer behaviour \ef{94} and hence, by matching with the
inner region asymptotics, lead to the asymptotic ODE \ef{124}. In
this class, the proof of Theorem \ref{Th.kap14} is
straightforward.

However, we must admit that we do not have a constructive way of
describing generic solutions. In fact, this is not that exciting
and/or surprising, since, even in the second-order case, solutions
of constant sign were attributed to generic ones {\sc only} by
using the Maximum Principle, which is not available for
bi-harmonic operators. For both second- and fourth-order parabolic
equations, conditions of attributing solutions of changing sign to
generic patterns are not fully known.

\ssk

\noi\underline{\sc Linear bi-harmonic equation: $\kappa=0$}.
However, the integrals in \ef{124} are, in general, {\em
oscillatory}, so that a proper regularity analysis becomes not
straightforward even in the linear case $\kappa=0$; see
\cite{GalPet2m}. Then the following holds:
 \be
 \label{13}
  \fbox{$
  \tex{
\kappa=0: \quad \int\limits^{\iy} G_4(\var(\t),0)\, {\mathrm d} \t
\quad \mbox{diverges to $-\iy$} \quad \Longleftrightarrow \quad
  a_0(\t) \to 0 \asA \t \to +\iy.
 }
 $}
  \ee
  Using asymptotic expansions of the  kernel \ef{i5} and the corresponding
   eigenfunctions, as well as sharp
values of the parameters \ef{i6},  yields a
  more practical condition:
  \be
  \label{14}
   \tex{
\frac {a_0'}{a_0} =  
   \var^{\frac 23}(\t) C_3\cos\big(b_0 \var^{\frac 43}(\t)+C_4\big)
   {\mathrm e}^{-d_0 \var^{4/3}(\t)}+...  \forA \t \gg 1,
    }
    \ee
 with some 
  constants $C_{3,4}$ depending in an obvious way on $C_{1,2}$ in \ef{i5}
and other parameters from \ef{i6}. Integrating yields
 \be
  \label{14NN}
   \tex{
\ln|{a_0(\t)}| = 
 \int\limits^\t
   \var^{\frac 23}(s) C_1 \cos\big(b_0 \var^{\frac 43}(s)+C_2\big)
   {\mathrm e}^{-d_0 \var^{4/3}(s)} \,{\mathrm d}s +...  \forA \t \gg
   1.
    }
    \ee
 The regularity condition \ef{13} is then re-formulated according to \ef{14NN}.
 Namely, the ``critical" backward parabola  occurs for the
 function (see \cite[\S~7]{GalPet2m})
  \be
  \label{var41}
 \var_*(\t) = 3^{-\frac 34}\, 2^{\frac {11}4}\, \big(\ln
\t\big)^{\frac
 34}+...
  \forA \t \gg 1,
 \ee
though, to guarantee divergence to minus infinity in \ef{13}, a
special ``oscillatory cut-off" of the function $\var_*(\t)$ may be
necessary.

 \ssk

\noi \underline{\sc Semilinear  equations}.
 For $\kappa \not =0$, instead of the linear  \ef{14}, we  deal
 with a
  nonlinear ODE
 \be
  \label{14N}
   \tex{
\frac {a_0'}{a_0} =  \hat \g
   \var^{\frac 23}(\t) C_3\cos\big(b_0 \var^{\frac
   43}(\t)+C_4\big)\,
   {\mathrm e}^{-d_0 \var^{4/3}(\t)}+ \kappa(a_0)+...  \forA \t \gg 1,
    }
    \ee
 and the analysis becomes more difficult.
 However, some of the results from Section \ref{S4.6} can be
 extended.
 Firstly, Proposition like \ref{Pr.99} can be restored provided
 an oscillatory cut-off of $\var(\t)$ is performed for the first
 integral in the right-hand side of \ef{14N} to be non-positive
 (though this business could look too artificial).
 Secondly,
 a statement similar to  Proposition \ref{Pr.reg2} remains valid
 with the ``linear" function \ef{a11} replaced by
 \be
 \label{a114}
 \tex{
\hat a_0(\t)=a_0(0) \, \eee^{\hat \g  C_3 \int_0^\t
   \var^{\frac 23}(s)\cos\big(b_0 \var^{\frac 43}(s)+C_4\big)
   {\mathrm e}^{-d_0 \var^{4/3}(s)}\, {\mathrm d} s},
   }
 \ee
 and with the corresponding changes in the integrals in \ef{a11},
 \ef{a12}.





   \ssk

 Let us briefly (and more formally) derive
 ``critical" nonlinearities $\kappa$. It is easy to see that
a somehow optimal and close to the critical dependence \ef{var41}
is then achieved for the nonlinearity
\be
\label{ll3}
 \tex{
  \kappa(v)=- \frac 1{|\ln v|} < 0 \forA v \approx 0^+.
  }
  \ee
 Indeed, solving the corresponding ODE without the linear term
 yields
 \be
 \label{ll1}
  \tex{
 \tilde a_0'= \kappa (\tilde a_0) \LongA
   \int_{\tilde a_0(\t)}^1 \frac{{\mathrm d} z}{z
   |\kappa(z)|}=\t \whereA \tilde
   a_0(\t) \to 0^+ \asA \t \to \iy.
   }
   \ee
  It follows from \ef{ll1} that, for the nonlinear coefficient
  \ef{ll3},
   \be
   \label{ll4}
   \tex{
   \tilde a_0(\t) ={\mathrm e}^{-\sqrt{2 \t}},
   }
   \ee
  so that, as is easy to see, that the linear term in negligible on the asymptotics \ef{ll4}, i.e.,
\be
   \label{ll2New}
    \tex{
    \var^{\frac 23}(\t)\,\, {\mathrm e}^{-d_0 \var^{4/3}(\t)}
    =o\big(|\kappa(\tilde a_0(\t))|\big) \asA \t \to \iy,
    \quad \mbox{provided that}
    }
    \ee
   \be
   \label{ll5}
    \var(\t) \gg (\ln \t)^{\frac 34} \forA \t \gg 1 \quad
    (\mbox{cf. (\ref{var41})}).
    \ee
 We thus arrive at a conclusion, which is similar to that in
 Proposition \ref{Pr.99}: for such negative $\kappa$'s, the vertexes of  arbitrarily ``wide"
 backward parabolae $\pa Q_0$ are regular.

Nevertheless, there are some principal differences with the much
simpler second-order case.
  For instance, if
 the integral in \ef{14NN} diverges and
 both linear and nonlinear terms on the right-hand side of
 \ef{14N} are sufficiently ``balanced", i.e., both equally
 involved in the asymptotics of $a_0(\t)$, the {\em actual checking
 regularity/irregularity of the origin becomes a principally non-solvable
problem}. It is curious that the most interesting ``interactional
case" (of linear and nonlinear terms in \ef{14N}) also begins at
functions such as \ef{kap99}, where the explicit constant $3 \sqrt
\pi$ must be replaced by a more complicated one composed from
those in \ef{i6} and $\g_{1,2}$  in \ef{114} uniquely given by the
BL-profile  \ef{z124}.

On the other hand, if the nonlinear term is asymptotically
negligible on the ``linear solutions" of \ef{14N}, then the
regularity and/or irregularity conditions remain practically the
same as for the pure bi-harmonic flow. These are rather trivial
results, which we do not intend to state and avoid such artificial
``rigorous" theorems.

 \subsection{Gradient dependent nonlinear perturbation}

For the second equation in \ef{sem4}, the rescaled equation in
\ef{ph44} takes the form
 \be
 \label{qq1}
  v_\t= \BB^* v+ \kappa(v) v \, (v_y)^4,
   \ee
so that using the same BL-profile \ef{124} and the variables
\ef{z24}, we obtain a similar dynamical system as in \ef{14N},
where the weaker nonlinear perturbation is estimates as in Section
\ref{S4.7}; cf. \ef{ss3}, where $\psi_0(y) \equiv F(y)$ is the
oscillatory kernel \ef{i5}. One can complete these computations;
however, as before, such gradient dependent nonlinear terms do not
affect the linear regularity criterion.

\subsection{Backward paraboloid in $\ren$}

Again, in greater detail, regularity analysis in $\ren$ for
Burnett equations (with the bi-Laplacian rather than the pure
Laplacian in the Navier--Stokes equations) is performed in
\cite[App.~A]{GalMazNSE}, so we present here a brief notice only.
For the equations  \ef{sem4N}, the lateral boundary of the domain
$Q_0$ in $\re^{N+1}$ can be given by the corresponding {\em
backward paraboloid}
 \be
 \label{parab14}
  \tex{
 \big( \sum_{i=1}^N a_i |x_i|^{2m} \big)^{\frac 1{2m}}= (-t)^{\frac 1{2m}} \,\, \var(\t), \quad \t= -
  \ln (-t), \quad a_i>0, \,\,\, \sum a_i^{2m}=1.
   }
   \ee
Again, a boundary layer study close to the rescaled (via \ef{z1})
boundary
 \be
 \label{vvv14}
  \tex{ \pa \hat Q_0: \quad
  \sum a_i |z_i|^{2m}=1,
  }
  \ee
  leads to a linear elliptic problem, which in the orthogonal
direction becomes ``quasi" one-dimensional, so that $g_0(\xi)$
given in \ef{z124} depends on the single variable \ef{norm1}.
Eventually, in Inner Region, the BL-behaviour leads to the
stabilization to
 \be
 \label{norm24}
 g_0(y,\t)=g_0\big(\var^{\frac 13}(\t)\, {\rm dist} \, \{y, \pa \tilde
 Q_0(\t)\}\big),
  \ee
  where $g_0(\xi)$ is as in \ef{z124}.
 This makes it possible to derive the asymptotic dynamical system
 for the first Fourier coefficient and hence an ODE regularity criterion for generic solutions.

  The resulting asymptotic ODE for $a_0(\t)$ is similar to \ef{124}, with
  the extra multiplier $\var^{N-1}$ in the first two terms on the
  right-hand side. Inevitably, the final ODE will depend
on the geometry of the backward paraboloid
  \ef{parab14} in a neighbourhood of its characteristic vertex
  $(0,0)$, which, in the most sensitive critical cases,
   makes it even less suitable for a definite regularity
  conclusion.

 \subsection{More on generalizations}

Using the above approach,
 there is no much principle differences and difficulties to treat the asymptotics
 of
 characteristic points for $2m$th-order poly-harmonic equations
  \be
  \label{2m11}
  u_t=(-1)^{m+1} \D^m u + f(x,t,u, \n u, D^2 u,...) \inB Q_0,
   \ee
   with zero Dirichlet (or others homogeneous) boundary conditions on $\pa Q_0$.
Though, of course, some involved technicalities occur indeed.
 Since the first rescaled variables are
  \be
  \label{sc66}
   \tex{
   u(x,t)=v(y,\t), \quad y= \frac x{(-t)^{1/2m}}, \quad \t = -
   \ln(-t),
   }
   \ee
most interesting  nonlinear terms in \ef{2m11} are now:
 $$
  \tex{
 f(\cdot)= \frac{1}{(-t)}\, \kappa(u)u, \quad \kappa(u)u \, |\n u|^{2m}, \quad
 \kappa(u) u |D^2 u|^m, \quad \mbox{etc.}
  }
 $$
 Then scalings \ef{sc66} lead to the rescaled parabolic equations
  \be
  \label{sc67}
  v_\t= \BB^* v + \left\{
   \begin{matrix}
    \kappa(v)v, \ssk\\
      \kappa(v)v \, |\n v|^{2m}, \ssk \\
 \kappa(v) v |D^2 v|^m,
 \end{matrix}
  \right. \quad \mbox{etc.,}
  \ee
   where $\BB^*$ is the linear adjoined operator \ef{B1*}. The
   corresponding dynamical systems for the expansion coefficients
   are obtained, as above, by (i) constructing a BL, and (ii) projecting the resulting
    PDEs \ef{sc67}, with the BL-approximation,  onto generalized Hermite polynomials
    \ef{psi**1}. The dynamical system, in general, becomes
    extremely oscillatory, and
     both linear and nonlinear terms can essentially affect regularity of the vertex
    $(0,0)$.

  \ssk

Finally, we again mention that here our main goal: to show how the
evolution of the first Fourier coefficient of generic solutions of
bi-harmonic PDEs leads to an ODE regularity criterion, has been
  We must admit however that, in some cases, this
did not   end up with constructive/deterministic regularity
conclusions, which are not always  possible and are even illusive
in general for higher-order nonlinear parabolic PDEs.

\end{small}
\end{appendix}

\end{document}